\documentclass[a4paper,11pt]{article}
\usepackage[hidelinks]{hyperref}
\hypersetup{
    colorlinks=true,
    linktoc=all,
    linkcolor=black,     
   citecolor=blue,     
}

\usepackage{tikz}
\usepackage{subfigure}
\usepackage{graphicx}
\usepackage{amssymb}
\usepackage{latexsym, bm}
\usepackage{multicol}
\usepackage{indentfirst}
\usepackage{amssymb,amsfonts}
\usepackage{amsmath}
\usepackage{hyperref}
\usepackage{setspace}
\newtheorem{theorem}{Theorem}[section]
\newtheorem{lemma}[theorem]{Lemma}
\newtheorem{claim}[theorem]{Claim}
\newtheorem{fact}[theorem]{Fact}

\newtheorem{proposition}[theorem]{Proposition}
\newtheorem{definition}[theorem]{Definition}
\newtheorem{problem}[theorem]{Problem}

\newtheorem{construction}[theorem]{Construction}
\newtheorem{conjecture}[theorem]{Conjecture}

\newtheorem{observation}[theorem]{Observation}
\newcommand{\ignore}[1]{}

\begin{document}

\begin{spacing}{0.97}
\date{}

\title{Two Ramsey-Tur\'{a}n numbers involving triangles\footnote{Center for Discrete Mathematics, Fuzhou University, Fuzhou, 350108, P.~R.~China. Correspondence Email: {\tt linqizhong@fzu.edu.cn}. Supported in part by NSFC (No.\ 12171088, 12226401) and NSFFJ (No. 2022J02018).}}

\author{Xinyu Hu \;\; \text{and} \;\; Qizhong Lin}
\maketitle

\begin{abstract}
Given integers $p, q\ge2$, we say that a graph $G$ is $(K_p,K_q)$-free if there exists a red/blue edge coloring of $G$ such that it contains neither a red $K_p$ nor a blue $K_q$. Fix a function $f( n )$, the Ramsey-Tur\'{a}n number $RT( {n,p,q,f( n ))} $ is the maximum number of edges in an $n$-vertex $(K_p,K_q)$-free graph with independence number at most $f( n )$.
For any $\delta>0$, let $\rho (p, q,\delta ) = \mathop {\lim }\limits_{n \to \infty } \frac{RT(n,p, q,\delta n)}{n^2}$. We always call $\rho (p, q):= \mathop {\lim }\limits_{\delta  \to 0}\rho (p, q,\delta )$ the Ramsey-Tur\'{a}n density of $K_p$ and $K_q$.

In 1993, Erd\H{o}s, Hajnal, Simonovits, S\'{o}s and Szemer\'{e}di proposed to determine the value of $\rho(3,q)$ for $q\ge3$, and they conjectured that for  $q \ge 2$, $\rho \left( {3,2q - 1} \right) = \frac{1}{2}(1 - \frac{1}{r(3,q) - 1})$. Recently, Kim, Kim and Liu (2019) conjectured that for  $q \ge 2$, $\rho( {3,2q } ) = \frac{1}{2}( 1 - \frac{1}{r( {3,q} )})$.
Erd\H{o}s et al. (1993) determined $\rho(3,q)$ for $q=3,4,5$ and $\rho(4,4)$. There is no progress on the Ramsey-Tur\'{a}n density $\rho (p, q)$ in the past thirty years.
In this paper, we obtain $\rho(3,6)=\frac{5}{12}$ and $\rho(3,7)=\frac{7}{16}$.
Moreover, we  show that the corresponding asymptotically extremal structures are weakly stable, which answers a problem of Erd\H{o}s et al. (1993) for the two cases.


\end{abstract}

\section{Introduction}

Ramsey theorem \cite{ram} states that for any integers $p_1,\dots,p_k$, there is a minimum integer, called Ramsey number $r=r(p_1,\dots,p_k)$, such that any $k$-coloring of edges of the complete graph $K_r$ contains a $K_{p_i}$ in the $i$th color for some $1\le i\le k$.
Subsequently, Tur\'{a}n posed a problem to determine the maximum number of edges of a $K_{p+1}$-free graph. In particular, Tur\'{a}n \cite{B9,turan54} proved that  the balanced complete $p$-partite graph, known as {\em Tur\'{a}n graph} $T_{n,p}$, is the unique extremal graph which attains the maximum number of edges among all $n$-vertex $K_{p+1}$-free graphs.
These Tur\'{a}n graphs have large independent sets of size linear in $n$, so it is natural to ask for the maximum number of edges of an $n$-vertex $K_{p+1}$-free graph without large independent set. Erd\H{o}s and S\'{o}s \cite{B6} initiated to study such Ramsey-Tur\'{a}n type problems, which have attracted a great deal of attention. We refer the reader to Simonovits and S\'os \cite{ss} for a nice survey on this topic, and \cite{bhs,bl,bls,bd,flz,fps,kkl,lr,sud} etc. for many recent results.

 Given integers $p_1,\dots,p_k$, we say that a graph $G$ is $(K_{p_1},\dots,K_{p_k})$-free if there exists a $k$-edge coloring of $G$ with no monochromatic copy of $K_{p_i}$ in the $i$th color for $1\le i\le k$. The {\em multicolor Ramsey-Tur\'{a}n number} $RT(n,p_1,\dots,p_k,m)$ is defined as the maximum number of edges of an $n$-vertex $(K_{p_1},\dots,K_{p_k})$-free graph $G$ with independence number $\alpha(G)<m$.
Clearly, there is no graph $G$ of order $n$ which is  $(K_{p_1},\dots,K_{p_k})$-free and $\alpha(G)<m$ if $n\ge r(p_1,\dots,p_k,m)$.

\begin{definition}[Ramsey-Tur\'{a}n density]\label{def1}
Given integers $p_1,\dots,p_k$ and $0<\delta<1$, let
\[
\rho (p_1,\dots,p_k,\delta ) = \mathop {\lim }\limits_{n \to \infty } \frac{RT(n,p_1,\dots,p_k,\delta n)}{n^2},
\]
 and $\rho (p_1,\dots,p_k): = \mathop {\lim }\limits_{\delta  \to 0} \rho (p_1,\dots,p_k,\delta )$ is the Ramsey-Tur\'{a}n density of graphs $K_{p_1},\dots,K_{p_k}$.

\end{definition}


For the case when $k=1$, the Ramsey-Tur\'{a}n density is well understood. For odd cliques, Erd\H{o}s and S\'{o}s \cite{B6} proved that $\rho( {{{2p + 1}}}) = \frac{1}{2}(1 - \frac{1}{p})$ for all $p \ge 1$. The problem for even cliques is much harder apart from the trivial case $K_2$. Erd\H{o}s and S\'{o}s \cite{B6} showed that $\rho ( {{4}}) \le \frac{1}{6}$. As an early application of the regularity lemma,  Szemer\'{e}di \cite{B8} showed that $\rho ( {{4}}) \le \frac{1}{8}$. No non-trivial lower bound on $\rho( {{4}} )$ was known until Bollob\'{a}s and Erd\H{o}s \cite{B1} provided a matching lower bound using an ingenious geometric construction, showing that $\rho ( {{4}}) = \frac{1}{8}$. Finally, Erd\H{o}s, Hajnal, S\'{o}s and Szemer\'{e}di \cite{B4} proved $\rho ( {{{2p}}} ) = \frac{1}{2}( 1 - \frac{3}{3p - 2})$ for all $p \ge 2$.

In general, it would be much more difficult to determine the exact values of $\rho(p, q)$ for $p,q\ge3$.
In 1993, Erd\H{o}s, Hajnal, Simonovits, S\'{o}s and Szemer\'{e}di \cite{B3} proposed the following problem.
\begin{problem}[Erd\H{o}s et al. \cite{B3}]\label{ehsi-p}
 Determine $\rho(3,q)$ for $q\ge3$.
\end{problem}

In particular, they made the following conjecture.

\begin{conjecture}[Erd\H{o}s et al. \cite{B3}]\label{ehsi-cj}
For $q \ge 2$, $\rho \left( {3,2q - 1} \right) = \frac{1}{2}(1 - \frac{1}{r(3,q) - 1})$.
\end{conjecture}

Recently, Kim, Kim and Liu \cite{kkl} obtained $RT(n,3,2q,g(n))=\frac{1}{2}(1 - \frac{1}{r(3,q)})n^2+o(n^2)$, where $g(n)={n}/{e^{\omega(n)\cdot(\log n)^{1-1/q}}}$ with $\omega(n)\to \infty$, and they also proposed the following conjecture.

\begin{conjecture}[Kim et al. \cite{kkl}]\label{kkl-cj}
For $q \ge 2$, $\rho \left( {3,2q} \right) = \frac{1}{2}(1 - \frac{1}{r(3,q)})$.
\end{conjecture}

However, we only know the following four classical Ramsey-Tur\'{a}n densities up to now.

\begin{theorem}[Erd\H{o}s et al. \cite{B3}]\label{01} We have
$\rho ( {{3},{3}} ) = \frac{1}{4}$, $\rho ( {{3},{4}} ) = \frac{1}{2}({1 - \frac{1}{3}})$, $\rho ( {{3},{5}} ) = \frac{1}{2}( {1 - \frac{1}{5}}),$  and $\rho ( {{4},{4}} ) = \frac{1}{2}( {1 - \frac{3}{14}}).
$
\end{theorem}


Given integers $p_1,\dots,p_k$, and the function $f(n)$, a sequence of graphs $\{G_n\}$ will be called an {\em asymptotically extremal} sequence of $RT(n,p_1,\dots,p_k,f( n ))$ provided $G_n$ is $(K_{p_1},\dots,K_{p_k})$-free, $\alpha ( {{G_n}} ) \le f( n )$, and
$e( {{G_n}} ) = RT( {n,p_1,\dots,p_k,f( n )} ) + o( {{n^2}} ).$

Given two graphs $S_n$ and $Z_n$, their distance $\Delta ( {{S_n},{Z_n}} )$  is defined as the minimum number of edges one has to delete from and add to $S_n$ to get a graph isomorphic to $Z_n$.
We shall say that for a problem $RT( {n,p_1,\dots,p_k,o(n)} )$ the {\em weak stability property} holds if for any two asymptotically extremal sequences $\{S_n\}$ and $\{Z_n\}$ of $RT( {n,p_1,\dots,p_k,o(n)} )$, their distance $\Delta ( {{S_n},{Z_n}} ) = o( {{n^2}} )$.

\begin{problem}[Erd\H{o}s et al. \cite{B3}]\label{Ed-st}
  Is it true that for $RT(n, {p_1}, \ldots ,{p_k},o(n))$ the asymptotically extremal structure is weakly stable?
\end{problem}

We know from \cite{B4} that for each $p\ge3$, the corresponding asymptotically extremal structure for $RT(n,p,o(n))$ is weakly stable. For the two color case, we know from \cite{B3} that for $RT(n,3,q,o(n))$ when $q=3,4,5$ and for $RT(n,4,4,o(n))$, the corresponding asymptotically extremal structures are weakly stable.

In this paper, we will determine $\rho(3,6)$ and $\rho(3,7)$, and show that the asymptotically extremal structures are weakly stable. This confirms Conjecture \ref{kkl-cj} and Conjecture \ref{ehsi-cj} for these two cases, and also resolves  Problem \ref{Ed-st} for the two cases.

\begin{theorem}\label{zhu}
We have $RT(n,3,6,o(n))=\frac{1}{2}(1 - \frac{1}{r(3,3)})n^2+o(n^2)$, i.e.,
$\rho ( {{3},{6}} ) = \frac{1}{2}(1 - \frac{1}{6}).$
Moreover, the asymptotically extremal structure for $RT(n,3,6,o(n))$ is weakly stable.
\end{theorem}

\begin{theorem}\label{zhu-2}
We have $RT(n,3,7,o(n))=\frac{1}{2}(1 - \frac{1}{r(3,4)-1})n^2+o(n^2)$, i.e.,
$\rho ( {{3},{7}} ) = \frac{1}{2}(1 - \frac{1}{8}).$
Moreover, the asymptotically extremal structure for $RT(n,3,7,o(n))$ is weakly stable.
\end{theorem}


\medskip

\noindent
{\bf  Notation:} For a graph $G=(V,E)$ with vertex set $V$ and edge set $E$, $e(G)$ is the number of edges  $|E|$ in $G$, and $uv$ or $(u,v)$ denotes an edge of $G$. For $X \subseteq V$, $e(X)$ is the number of edges in $X$, and $G[X]$ denotes the subgraph of $G$ induced by $X$. For disjoint $X_1,\dots,X_t\subset V$, let $G[X_1,\dots,X_t]$ denote the subgraph induced by all edges between them. For two disjoint subsets $X,Y\subseteq V$, $e_G(X,Y)$ denotes the number of edges between $X$ and $Y$, and $N_G(X,Y)$ denotes the common neighborhood of vertices of $X$ in $Y$, and let $\deg_G(X,Y)=|N_G(X,Y)|$. If $X=\{v\}$, we will use $N_G(v,Y)$ and $\deg_G(v,Y)$ instead of $N_G(X,Y)$ and $\deg_G(X,Y)$, respectively. In particular, the neighborhood of a vertex $v$ in $G$ is denoted by $N_G(v)=N_G(v,V)$ and the degree of $v$ in $G$ is $\deg_G(v)=|N_G(v)|$.
Let $X \sqcup Y$ denote the disjoint union of $X$ and $Y$. A complete $p$-partite graph with vertex set $\sqcup_{i=1}^pV_i$, where $|V_i|=n_i$, is denoted by $K_p(n_1, \dots, n_p)$. Let $[n]=\{1,2,\dots,n\}$, and $[m,n]=\{m,m+1,\dots,n\}$. We always delete the subscriptions if there is no confusion from the context.

\medskip\noindent
{\bf Organization:} In Section \ref{pre}, we will include many essential constructions of graphs together with useful lemmas.
 In Section \ref{main1} and Section \ref{k3-k7}, we will present the proofs of Theorem \ref{zhu} and Theorem \ref{zhu-2}, respectively. Finally,
we will mention some interesting problems in Section \ref{crp}.

\section{Preliminaries}\label{pre}

\subsection{Szemer\'{e}di regularity lemma}
In 1975, Szemer\'{e}di \cite{sze75} confirmed the famous conjecture of Erd\H{o}s and Tur\'{a}n \cite{et} stating that any sequence of integers with positive upper density must contain arithmetic progressions of any finite
length, in which Szemer\'{e}di used a deep and complicated combinatorial lemma. The lemma now is always called Szemer\'{e}di regularity lemma, which is a powerful tool in extremal graph theory.
 There are many important applications of the regularity lemma, we refer the reader to nice surveys \cite{kss,ko-sim,rs} and other related references.

Let $X,Y\subseteq V(G)$ be disjoint nonempty sets of vertices in a graph $G$. The density of $(X,Y)$ is $d_G(X,Y)=\frac{e_G(X,Y)}{|X||Y|}$. For $\varepsilon>0$, the pair $(X,Y)$ is $\varepsilon$-regular in $G$ if for every pair of subsets $X'\subseteq X$ and $Y'\subseteq Y$ with $|X'|\geq \varepsilon |X|$ and $|Y'|\geq \varepsilon |Y|$ we have $| {{d_G}( {X,Y} ) - {d_G}( {X',Y'} )} | \le \varepsilon $. Additionally, we say that $(X,Y)$ is $(\varepsilon,\gamma)$-regular if $d_G(X,Y)\geq\gamma$ for some $\gamma >0$. We say a partition $V(G)=\sqcup_{i=0}^m V_i$ of $G$ is equitable if $|V_0|\le\varepsilon n$ and $||V_i|-|V_j||\le 1$ for all distinct $i$ and $j$ in $[m]$. An equitable partition $V(G)= \sqcup_{i=0}^m{{V_i} } $
is an $\varepsilon$-regular partition of a $k$-edge-colored graph $G$ if for each $i\in[m]$, all but at most $\varepsilon m$ choices of $j\in[m]$ satisfy that the pair $(V_i,V_j)$ is $\varepsilon$-regular in $G_\ell$ for each color $\ell\in[k]$.

In order to show the weak stability properties, we will use the following regularity lemma.
\begin{lemma}[Szemer\'{e}di \cite{sze78}]\label{reg-le}
Suppose $0<\frac{1}{M'}\ll \varepsilon$, $\frac{1}{M}\ll\frac{1}{k}\leq1$, and $n\geq M$. Suppose that $G$ is an $n$-vertex $k$-edge-colored graph and $U_1\sqcup U_2$ is a partition of $V(G)$. Then there exists an $\varepsilon$-regular equitable partition $V(G)= \sqcup_{i=0}^m{{V_i} } $ with $M\leq m\leq M'$ and for each $i \in[m]$, we have either $V_i \subseteq U_1$ or $V_i \subseteq U_2$.
\end{lemma}

We have two basic properties of the regularity pair as follows.
\begin{lemma}\label{M-s-d}
Let $(X,Y)$ be an $\varepsilon$-regular pair of density $d$ in a graph $G$. If $B\subseteq Y$ with $|B|\ge \varepsilon|Y|$, then there exists a subset $A\subseteq X$ with $|A|\ge (1-\varepsilon)|X|$ such that each vertex in $A$ is adjacent to at least $(d-\varepsilon)|B|$ vertices in $B$.
\end{lemma}

\begin{lemma}\label{Sli-L}
Let $\varepsilon<\alpha,\gamma,\frac{1}{2}$. Suppose that $(X,Y)$ is an $(\varepsilon,\gamma)$-regular pair in a graph $G$. If $X'\subseteq X$ and $Y'\subseteq Y$ satisfy $|X'|\ge \alpha|X|$ and $|Y'|\ge \alpha|Y|$, then $(X',Y')$ is an $(\varepsilon',\gamma-\varepsilon)$-regular pair in $G$, where $\varepsilon':=\max\{\varepsilon/\alpha,2\varepsilon\}$.
\end{lemma}

We will also use the following lemma by Balogh,  Liu and Sharifzadeh \cite[Lemma 3.1]{bls} which refines a result of Erd\H{o}s et al. \cite[Lemma 2]{B3}. We include the short proof for completeness.
\begin{lemma}[Balogh, Liu and Sharifzadeh \cite{bls}]\label{b-indent}
Let $G$ be an $n$-vertex graph with $\alpha(G)\leq \delta n$ for some $0<\delta<1$, and let $\varphi:E(G)\rightarrow [2]$. Then there is a partition $V(G)=U_1\sqcup U_2$ such that for every $k\in [2]$, $\alpha(G_k[U_k])\leq \sqrt{\delta}n.$
\end{lemma}
{\em Proof.} Let $X_0 = V(G)$ and $Y_0 = \emptyset$. At step $i\ge0$, if $\alpha(G_1[X_{i}])\le \sqrt{\delta}n$, then we stop. Otherwise, there is a subset $S_{i}\subseteq X_{i}$ with $|S_i|>\sqrt{\delta}n$ and $S_i$ contains no edge in color 1. Together with $\alpha(G)\leq \delta n$, we have $\alpha(G_2[S_i])\leq \delta n$.
Define $X_{i+1}:= X_i\setminus S_i$ and $Y_{i+1}:= Y_i\cup S_i$.
Clearly, $\alpha(G_2[Y_{i+1}])\leq \alpha(G_2[Y_{i}])+\delta n$. Suppose the iteration stops after $t$ steps, then $t\le1/\sqrt{\delta}$ as we delete at least $\sqrt{\delta}n$ vertices each time from $X_i$. Finally, we have $\alpha(G_1[X_{t}])\le \sqrt{\delta}n$, and $\alpha(G_2[Y_{t}])\leq t\cdot\delta n \le \sqrt{\delta}n$, and so $X_t$ and $Y_t$ form a partition of $V(G)$ as desired.\hfill$\Box$

\subsection{Useful constructions and weighted graphs}\label{petg}
In the following, we will introduce many related constructions of graphs and useful lemmas.
Let us begin with a geometric construction by Erd\H{o}s and Rogers \cite{B5}.

\medskip\noindent
\textbf{Erd\H{o}s graph} (or called Erd\H{o}s-Rogers graph):  There are a constant $c>0$ and $n_0$ such that for every $n>n_0$ there exists an $n$-vertex graph $G_n$ satisfying
$K_3 \nsubseteq {G_n}$, and $\alpha (G_n) \le n^{1-c}.$

\smallskip
{\em Remark.} 
In \cite{alon1},  Alon constructed an $n$-vertex graph $G_n$ that is $K_3$-free and $\alpha (G_n) =O(n^{2/3}).$ For more constructions, see \cite{alon2,kpr} and the related references therein.
Indeed, we can take $G_n$ such that it is $K_3$-free and $\alpha (G_n) = O(\sqrt{n\log n})$ from the celebrated result of Kim \cite{kim}. 

\smallskip

 Given a graph ${G}$ whose vertex set  is partitioned into the classes $X_1, \ldots ,{X_q}$, an $2$-coloring of edges of $G$ will be called a {\bf canonical coloring} if all the edges lie in ${X_i}$ have the same color; moreover, all the edges between distinct $X_i$ and $X_j$ have the same color.

\begin{construction}[Erd\H{o}s et al. \cite{B3}]\label{con1}
Let $t = r(p,q) - 1$. Color the Tur\'{a}n graph ${T_{n,t}}$ by red and blue canonically (with respect to the classes of ${T_{n,t}}$) so that the colored graph contains neither a  red $K_q$, nor a  blue $K_p$. Put into each class of this graph a {\bf red} Erd\H{o}s graph. Then the resulting colored graph $U(n,p,q)$ contains neither a  red $K_{2q-1}$, nor a  blue $K_p$. In particular, $\rho \left( {p,2q - 1} \right) \ge \frac{1}{2}(1 - \frac{1}{r(p,q) - 1}).$
\end{construction}

\ignore{From the above construction, we immediately have
\begin{align}\label{low-odd-q}
\rho \left( {p,2q - 1} \right) \ge \frac{1}{2}\left(1 - \frac{1}{r(p,q) - 1}\right).
\end{align}}

For the even case, we have the following lower bound.
\begin{lemma}\label{Lower bound}
  Let $t = r(p,q)$ for $p,q\ge3$. Then
 $
RT( {n,p,2q,o(n)} ) \ge\left ( {1 - \frac{1}{t}}\right ){n\choose2}+\frac{\rho(p)}{t^2}n^2+o(n^2).
$
Or equivalently,  $\rho(p,2q)\ge\frac12\left(1 - \frac{1}{t}\right)+\frac{\rho(p)}{t^2}.$
\end{lemma}
\noindent{\em Proof.} Let us construct a colored graph $H( {n,p,q} )$ as follows.
Let $G= {T_{n,t}}$ with $t$ classes $X_1, \ldots,X_t$. Since $t = r(p,q)$, we can color $G_1=G[ {{X_1}, \ldots ,{X_{t-1}}} ]$ by red and blue canonically (with respect to the classes $X_1, \ldots ,X_{t-1}$) so that the colored graph contains neither a  red $K_{q}$ nor a  blue $K_p$. Then we put into each class ${X_i}$ for $i\in[t-1]$ a {\bf red} Erd\H{o}s graph, so the resulting graph contains neither a  red $K_{2q-1}$ nor a blue $K_p$. Finally, we put into ${X_t}$ a {\bf blue} extremal graph for $RT(|X_t|,p,o(n))$ and color all edges between $\sqcup_{i=1}^{t-1}X_i$ and ${X_t}$ {\bf red}. The final colored graph, denoted by
${{{H_n} = H( {n,p,q} )}},$ contains neither a  red $K_{2q}$ nor a  blue $K_p$.
Clearly, we have $\alpha ( {{H_n}} )= o( n )$, and so
$
RT( {n,p,2q,o(n)} ) \ge \left ( {1 - \frac{1}{t}}\right ){n\choose2}+\rho(p)|X_t|^2+o(n^2)
$ as claimed.
\hfill$\Box$




\medskip
In \cite{B1}, Bollob\'{a}s and Erd\H{o}s gave an ingenious geometric construction, showing that $\rho ( {{4}} ) = \frac{1}{8}$. The constructed graph, called BE-graph, is as follows.

\medskip\noindent
{\bf BE-graph}: For any $\varepsilon>0$ and a large integer $h$, we fix a sufficiently large $n_0(\varepsilon, h)$ and assume that $n>n_0(\varepsilon, h)$ is even and $\mu=\frac{\varepsilon}{\sqrt{h}}$. Partition an $h$-dimensional unit sphere {\bf S} into $n/2$ domains $D_1,\dots,D_{n/2}$, of equal measure and diameter less than $\mu/10$. Let $X_1$ and $X_2$ be two disjoint sets of size $n/2$, each containing exactly one point from each $D_i$ for $i\in [\frac n2]$. The $BE$-graph is defined on $X_1\sqcup X_2$,  where $(x,y)$ is an edge iff one of the following holds:

(i) $x\in X_1$ and $y\in X_2$, and the distance $d(x,y)<\sqrt{2}-\mu$.

(ii) $x,y\in X_i$ for $i=1$ or $2$, and the distance $d(x,y)\ge2-\mu$.

\smallskip
We remark that the $BE$-graph is $K_4$-free with $\alpha(BE)=o(n)$ and each vertex has degree roughly $n/4$, and the subgraph spanned by $X_i$ for $i\in[2]$ is $K_3$-free with edge number $o(n^2)$.

In \cite{B3}, Erd\H{o}s et al. gave a generalization of the $BE$-graph as follows. One can see \cite{lrss} for a class of $BE$-type-graphs on the high dimensional {\em complex} sphere, which can be used to give all rational densities for generalized Ramsey-Tur\'{a}n densities.

\begin{construction}[Weighted $t$-partite {\em GBE}-graph \cite{B3}]\label{con2}
 We define a graph $$B:=B(h, t|n_1, \ldots ,n_t|\mu),$$
 which depends on the integers $h$, $t$, ${n_1}, \ldots ,{n_t}$, and a small number $\mu>0$.  We shall subdivide the points of the $h$-dimensional sphere \textbf{S} into $n_i$ sets of equal measure and diameter $ < \frac{\mu }{{10}}$, for each $i\in[t]$. For $i \in [t]$ and $\ell \in [n_i]$, we choose a vertex ${x_{i\ell}}$ from the $\ell$th set of the $i$th partition. Put ${X_i} = \{ {{x_{i\ell}}: \ell \in [n_i]} \}$ for $i \in [t]$ and let $V( B ) = \sqcup_{i = 1}^t {{X_i}} $.
 (1) For each pair $x,y \in {X_i}$ we join them iff $d( {x,y} ) > 2 - \mu $.
 (2) For some pairs $( {i,j} )$ $( {i \ne j} )$ we join every vertex of $X_i$ to every vertex of $X_j$,
 (3) for all the other pairs we build a $BE$-graph between the two classes: join $x \in {X_i}$ to $y \in {X_j}$ iff $d( {x,y} ) < \sqrt 2  - \mu $.
 The resulting graph is the generalized Bollob\'{a}s-Erd\H{o}s or shortly $GBE$-graph. 
\end{construction}



When finding the desired $K_p$, we will use the so-called {\em generalized complete graph}, in which each ``marked" vertex (with desired color) always contributes two vertices to it.

\begin{definition}[Generalized complete graph \cite{B3}]
\label{def4}
Let $H$ be a graph in which some vertices are ``marked", and the edges are assigned with weights $0$, $\frac{1}{2}$ or $1$. $X \subseteq V(H)$ and $Y \subseteq X$ span a generalized complete graph of the size $| X | + | Y |$ if
(a) all the vertices of $Y$ are ``marked",
(b) all the edges in $Y$ have weight $1$, and
(c) all the edges in $X$ have weight $ \ge \frac{1}{2}$.
\end{definition}


\begin{definition}[Weighted Ramsey numbers \cite{B3}]\label{wrn}
 Assume that the vertices and the edges of a $K_t$ are colored by the colors of $[2]$, and the edges are weighted by $\frac{1}{2}$ or $1$. For the color $k\in[2]$ we define ${\Lambda_k}$ as a weighted graph, where a vertex is ``marked" iff its color is $k$, and the edge $e$ gets weight $\frac{1}{2}$ or $1$ if its color is $k$ and its weight is $\frac{1}{2}$ or $1$ respectively.  Finally, assume that there is a distribution $ {\bf u}  = ( {{u_1}, \ldots ,{u_t}} ) \; ( {{u_i} \ge 0,\sum_{i = 1}^t {{u_i}}  = 1})$ on the vertices of $K_t$. If we wish to emphasize that $K_t$ is weighted, we shall write $K_t(\omega)$. Let
 $g( {{K_t}( \omega  ),\bf{u} } ) = \sum_{1 \le i < j \le t} {{\omega(i,j)}{u_i}{u_j}},$ where ${{\omega(i,j)}}$ is the weight of $( {i,j} )$. We define the edge density of such a weighted graph as
$$g(K_t( \omega )) =\max\limits_{\bf{ u} } g( {{K_t}( \omega),\bf{u} }) =\max\limits_{\bf{ u} } \sum\limits_{1 \le i < j \le t} {{\omega(i,j)}{u_i}{u_j}}.$$
The weighted Ramsey number $\beta  (p, q)$ is the maximum $\beta$ such that there exist an integer $t$ and a weighted coloring of $K_t( \omega )$ with edge density $\beta$ for which $\Lambda_1$ contains no generalized $K_p$ and $\Lambda_2$ contains no generalized $K_q$. Clearly, $t < r(p, q)$.
\end{definition}

\begin{lemma}[Erd\H{o}s et al. \cite{B3}]\label{not generalized}
The canonically colored  $GBE$-graph $B( {h, t |{n_1}, \ldots , {{n_t}} |\mu } )$ does not contain a monochromatic $K_\ell$ iff the corresponding colored weighted $K_t(\omega)$ does not contain a generalized $K_\ell$ in the same color.
\end{lemma}

We remark that a vertex of a generalized $K_\ell$ is ``marked" with color $k$ means the edges in the corresponding part of the corresponding {\em GBE}-graph are all colored $k$.

Let $A=(a_{ij})_{t\times t}$ be a $t\times t$ matrix, whose diagonal entries are all $0'$s, and $a_{ij}=\frac{1}{2}$ or $1$ depending on whether the $i$th and $j$th classes of {\bf GBE}-graph are joined by a $BE$-graph or completely. Indeed, $A$ is the weighted adjacent matrix of $K_t(\omega)$ corresponding to {\em GBE}-graph.

Let us normalize the integers $n_1, \ldots  , n_t$ by putting $ {u_i} = \frac{n_i}{n},$ where $n=\sum_{i=1}^tn_i$.
 Let \textbf{u} be a $t$-dimensional vector whose $i$th component is $u_i$.
Note that
\begin{align}\label{B-g}
e(GBE)=\frac{1}{2}\textbf{u}A{\textbf{u}^T}{n^2}= g(K_t(w),{\bf{u}})n^2+o(n^2)
\end{align}
as $h\rightarrow\infty ,\varepsilon \rightarrow 0, n\rightarrow \infty$.

If $\textbf{u}A{\textbf{u}^T}$ attains its maximum
on the boundary, i.e., for some $u_k=0$, then the $GBE$-graph will be called \textbf{degenerate}. If there are no maxima on the boundary, then one can easily see that there is exactly one maximum, and the structure of this $GBE$-graph will be called \textbf{dense structure}. Then the corresponding $K_t(\omega)$ is said to be {\bf dense} as well.


The following result gives the structure of the asymptotically extremal graphs.
\begin{theorem}[Erd\H{o}s et al. \cite{B3}]\label{00}
Given  $p, q\ge 3$, then for some fixed $t$ there exists a sequence of dense GBE-graphs $B( {h, t |{n_1}, \ldots , {{n_t}} |\mu } )$ asymptotically extremal for $RT( {n,p, q,o( n )} )$. In particular,
$
\rho (p,q) = \beta (p,q)=g(K_t(\omega)),
$ where $K_t(\omega)$ corresponding to GBE-graph.
\end{theorem}

We also need the following simple lemma.

\begin{lemma}\label{equal}
Let ${B_n}$ be a graph with $e( {{B_n}} ) = ( {1 - \frac{1}{{k + 1}}} ){n\choose2}+o(n^2)$, and $A\subseteq V(B_n)$ and let $m = n - | A |$. If $A$ contains $o( {{n^2}} )$ edges, and $e( {{B_n} - A} ) \le ( {1 - \frac{1}{k }} ){m\choose2}+o(n^2)$, then we have $m=\frac{k}{k+1}n+o(n)$ and $e( {{B_n} - A} ) = ( {1 - \frac{1}{k }} ){m\choose2}+o(n^2)$.
\end{lemma}
\noindent
{\em Proof.}~Since $e( {{B_n}} ) = e( {B_n} - A ) + e(A,V(B_n)\setminus A ) + e( A )$, it follows that
\begin{align*}
  e( {{B_n}} )  \le \left( {1 - \frac{1}{k}} \right){m\choose2} + m( {n - m} ) + o( {{n^2}} )
  \le\left( {1 - \frac{1}{{k + 1}}} \right){n\choose2}  + o\left( {{n^2}} \right),
  \end{align*}
  where the last inequality holds since $({1 - \frac{1}{k}} )\frac{{{x^2}}}{2} + x( {n - x} )$ takes maximum value at $x=\frac{k}{k+1}n$.
  Thus, the assertion follows since $e( {{B_n}} ) = ( {1 - \frac{1}{{k + 1}}} ){n\choose2}+o(n^2)$ from the assumption.\hfill$\Box$

\medskip
The proof of the above lemma already implies the following lemma.
\begin{lemma}[Erd\H{o}s et al. \cite{B3}]\label{10}
Let ${B_n}$ be a graph and $A\subseteq V(B_n)$, and let $m = n - | A |$. If $A$ contains $o( {{n^2}} )$ edges, and
$
e( {{B_n} - A} ) \le ( {1 - \frac{1}{k}} ){m\choose2},
$
then
$
e( {{B_n}}) \le ( {1 - \frac{1}{{k + 1}}}){n\choose2}.
$
Or equivalently,
if for a vertex $x \in {K_t}( \omega  ) $,
$
g( {{K_t}( \omega  ) - x} ) \le \frac{1}{2}( {1 - \frac{1}k} ),
$
then
$
g( {{K_t}( \omega  )} ) \le \frac{1}{2}( {1 - \frac{1}{k+1}}).
$
\end{lemma}

\medskip For a (colored) graph $G_t$ of order $t$ with the distribution \textbf{u} on the vertices and weights \textbf{w} on the edges, let the \textbf{weighted degree} of a vertex $x$ be
\begin{align}\label{wt-d}
d(x) = t\cdot\sum\limits_{(x,y)\in E( {{G_t}} )} {{\omega(x,y)}{u_y}},
\end{align}
 where $\omega(x,y)$ is the weight of an edge $(x,y)$ of $G_t$.

 \medskip
The following result describes the optimal weight distribution of vertices of a dense ${K_t}( \omega  )$.
\begin{lemma}[Brown, Erd\H{o}s and Simonovits \cite{B2}]\label{0}
The optimal weight distribution of vertices of a dense ${K_t}( \omega  )$ is attained when all the weighted degrees are equal. In particular, $g({K_t}(\omega))=d(x)/(2t)$, where $x$ is any vertex of ${K_t}( \omega  )$.
\end{lemma}
\noindent{\em Proof.} Note that the optimal weight distribution of vertices of a dense ${K_t}( \omega  )$ is attained in the internal domain, so the assertion follows from the Lagrangian method and the definition of the weighted degree immediately.\hfill$\Box$

\begin{lemma}[Erd\H{o}s et al. \cite{B3}]\label{123}
If in a dense $ {K_t}( \omega  ) $ each vertex is incident to at least $\lambda $ half  edges, then
$2g( {{K_t}( \omega  )} ) \le 1 - \frac{1}{t} - \frac{\lambda }{{2t}}.$
If there exists no $\lambda$-regular graph on $t$ vertices, then we have strict inequality $2g( {{K_t}( \omega  )} ) <1 - \frac{1}{t} - \frac{\lambda }{{2t}}.$
\end{lemma}

In ${K_t}(\omega)$, we say an edge is {\bf full} or {\bf half} if it has weight 1 or $1/2$, respectively. We say that a graph $H$ is {\bf full/half} if all edges of $H$ are full/half.
For a subset $X$ of vertices of ${K_t}(\omega)$, we say that a vertex $u$ of ${K_t}(\omega)$ is {\bf full/half-adjacent} to $X$ if all edges $xu$  for $x\in X$ are full/half. The {\bf full-neighborhood} of $u$ is the subset of ${K_t}(\omega)$ which is full-adjacent to $u$.


\section{Proof of Theorem \ref{zhu}}\label{main1}

\subsection{$\rho(3,6)=5/12$}\label{3.2}
The lower bound of $\rho(3,6)\ge5/12$ follows from Lemma \ref{Lower bound} immediately, and so it suffices to show $\rho(3,6)\le5/12$.


Let us first have the following lemma.
\begin{lemma}\label{1}
Let, in ${K_t}( \omega )$, $p$ disjoint ${C_5}$'s, $q$ disjoint ${K_4}$'s, $r$ disjoint ${K_3}$'s and $s$ independent edges be given, where these $s$ independent edges and that of ${C_5}$'s , ${K_4}$'s and ${K_3}$'s are also pairwise vertex disjoint and all are half. Then $2g( {{K_t}( \omega  )} ) \le 1 - \frac{{30}}{{30t - 75p - 72q - 45r - 20s}}.$
\end{lemma}
\noindent
{\em Proof.} Let ${K_t}( \omega )$ be defined on $[t]$. Note that the density will not decrease if we replace a half  edge by a full  edge,  so  we may assume that all the other edges are full  edges. There are $w = t - 5p - 4q - 3r - 2s$ ``other vertices". This structure will be denoted by $W( {p,q,r,s,w} )$.

\begin{claim}\label{dense}
Let $A=(a_{ij})_{\ell\times \ell}$ and $A'=(a_{ij}')_{\ell'\times \ell'}$ be two matrixes corresponding to two dense structures of GBE-graphs and let structure $C$ be obtained by joining each vertex of $A$ to each vertex of $A'$ with weight $1$, then the structure $C$ is also a dense structure.
\end{claim}
\noindent{\em Proof.} Let ${\bf{u}}=({\bf u_1 \ u_2})$, where ${\bf{u_1}}=(x_1, \ldots ,x_\ell)$ and ${\bf{u_2}}=(y_1, \ldots ,y_{\ell'})$ be $1\times \ell$ and $1\times \ell'$ vectors, respectively. Let $J$ be the $\ell'\times \ell$ all $1$ matrix. Note that
\[{\bf{u}}C{{\bf{u}}^T} = ({\bf u_1 \ u_2})\left({\begin{array}{*{20}{c}}
A&{{J^T}}\\
J&A'
\end{array}} \right)\left({\begin{array}{*{10}{c}}
{{\bf{u_1}}^T}\\
{{\bf{u_2}}^T}
\end{array}}\right) = {\bf{u_1}}A{\bf{u_1}}^T + {\bf{u_2}}A'{\bf{u_2}}^T + 2\Bigg( {\sum\limits_{1\le i\le\ell} {{x_i}} } \Bigg)\Bigg( {\sum\limits_{1\le j\le\ell'} {{y_j}} } \Bigg).\]

On the contrary, suppose that the structure $C$ attains the maximum at the boundary, say at ${\bf{u_0}}=(x_1,\dots,x_\ell,y_1,\dots,y_{\ell'})$ with $x_\ell=0$ and $\sum_{i = 1}^\ell {{x_i}}  + \sum_{j = 1}^{\ell'} {{y_j}}  = 1$. Let $a:=\sum_{i = 1}^\ell {{x_i}}$. We may assume that $A$ attains its maximum at $(x_1',\dots,x_{\ell}')$ with $x_i'\neq0$ for $i\in[\ell]$ and $\sum_{i = 1}^\ell {{x_i'}}=1$ since the structure $A$ is dense. Then ${\bf{u_1}}A{{\bf{u_1}}^T}$ will take larger value at $(ax_1',\dots,ax_\ell')$ than that at $(x_1,\dots,x_\ell)$. Together with $\sum_{i = 1}^\ell {a{x_i'}}=\sum_{i = 1}^\ell {{x_i}}=a$, we have that ${\bf{u}}C{{\bf{u}}^T}$ will take larger value at $(ax_1',\dots,ax_\ell',y_1,\ldots,y_{\ell'})$ than that at ${\bf{u_0}}$, which is impossible. The claim follows.   \hfill$\Box$

\smallskip
Let $A$ be the corresponding (weighted) adjacent matrix of the subgraph induced by the vertices of a half $C_5$.
\ignore{Then $$A=\left( {\begin{array}{*{20}{c}}
0&{1/2}&1&1&{1/2}\\
{1/2}&0&{1/2}&1&1\\
1&{1/2}&0&{1/2}&1\\
1&1&{1/2}&0&{1/2}\\
{1/2}&1&1&{1/2}&0
\end{array}} \right).$$}
Then
$$A=\left( {\begin{array}{*{20}{c}}
0&{\frac{1}{2}}&1&1&{\frac{1}{2}}\\
{\frac{1}{2}}&0&{\frac{1}{2}}&1&1\\
1&{\frac{1}{2}}&0&{\frac{1}{2}}&1\\
1&1&{\frac{1}{2}}&0&{\frac{1}{2}}\\
{\frac{1}{2}}&1&1&{\frac{1}{2}}&0
\end{array}} \right).$$
Let ${\bf{u}} = \left( {{x_1}, \ldots ,{x_5}} \right)$ (${x_i} \ge 0$ and $\sum_{i = 1}^5 {{x_i}}  = 1$), then
${\bf{u}}A{\bf{u}}^T$ takes the maximum $0.6$ when $x_i=0.2$ for all $i\in[5]$, so all half $C_5$'s as assumption are dense structures. Similarly, all these half $K_4$'s, $K_3$'s, $K_2$'s, and ``other vertices'' as assumption are dense structures.
Hence, $W:=W( {p,q,r,s,w} )$ is a dense structure from Claim \ref{dense}.

Further, suppose the optimum distribution of $W$ is $x, \ldots ,x$, $y, \ldots ,y$, $z, \ldots ,z$, $u, \ldots ,u$, $v, \ldots ,v$, then $5px+4qy+3rz+2su+wv=1$. For a $C_5$ on vertex set $[5]\subset [t]$,  note that the edges of this $C_5$ are half and all the other edges are full from the assumption, then from (\ref{wt-d}) for each $i\in[5]$, the weighted degree
\begin{align*}d(i) =t\left(\sum\limits_{k \in [ t ]\backslash [ 5 ]} {{w{( {i,k} )}}{u_k}}  +
  \sum\limits_{k \in [ 5 ]\backslash \{ i\} } {{w{( {i,k} )}}{u_k}}\right)
  = t\big[1 - 5 x + \left(2\cdot x/2 + 2x\right)\big]
   =t(1-2x),
   \end{align*}
\ignore{
\begin{align*}
  d( i ) =t\Bigg(\sum\limits_{k \in [ t ]\backslash [ 5 ]} {{w{( {i,k} )}}{u_k}}  +
  \sum\limits_{k \in [ 5 ]\backslash \{ i\} } {{w{( {i,k} )}}{u_k}}\Bigg)
  = t\left[1 - 5 x + \left(2\cdot x/2 + 2x\right)\right]
   =t(1-2x).
\end{align*}}
Similarly, the weighted degrees of the vertices of $K_4$, $K_3$, $K_2$, and
``other vertices'' are $t(1 - 5y/2)$, $t(1-2z)$, $t(1-3u/2)$, $t(1 - v)$, respectively. By Lemma \ref{0}, all the weighted degrees are equal.
 Hence, if we set $x = 15m$, then $y = 12m$, $z = 15m$, $u = 20m$, $v = 30m$. Therefore,
 \begin{align*}
   1 & = 5p\cdot15m + 4q\cdot12m + 3r\cdot15m + 2s\cdot20m + ( {t - 5p - 4q - 3r - 2s} )\cdot30m \\
    & = ( {30t - 75p - 72q - 45r - 20s} )m.
 \end{align*}
By Lemma \ref{0},
$2g({{K_t}( \omega  ))}  =1-2x=1-30m = {1 - \frac{{30}}{{30t - 75p - 72q - 45r - 20s}}}$ as desired.\hfill$\Box$

\medskip

We call a red/blue edge coloring of $K_5$ a {\bf pentagonlike coloring} if each of the subgraphs induced by all red or blue edges is $C_5$.
The following lemma is clear.
\begin{lemma}\label{jian}
Let $\phi $ be a red/blue edge coloring of $K_5$ such that it contains neither a red $K_3$ nor a blue $K_3$, then $\phi $ must be the pentagonlike coloring.
\end{lemma}

From Theorem \ref{00}, let $B_n(t)=B(h, t |{n_1}, \ldots , {{n_t}} |\mu)$, which is a dense asymptotically extremal graph of $RT( {n,{3},{6},o( n )} )$, with its {\bf canonical coloring} with respect to the $t$ classes $X_1,\dots,X_t$ of $B_n(t)$. Let ${K_t}( \omega  )$ be the corresponding weighted colored complete graph with distribution ${\bf{u}}=(u_1,\dots,u_t)$ \;$( {{u_i} > 0,\sum_{i = 1}^t {{u_i}}  = 1})$, in which each edge is {\bf full} or {\bf half}. We say a {\bf vertex} of ${K_t}( \omega  )$ is {\bf blue (red)} if all edges in the corresponding class $X_\ell$ for some $\ell\in[t]$ are blue (red).

From the assumption, the colored $B_n(t)$ contains neither a  blue $K_3$ nor a  red $K_6$ and $\alpha(B_n(t))=o(n)$. We may assume that $e(B_n(t))\ge\frac{1}{2}( {1 - \frac{1}{6}} ){n^2} + o(n^2)$ from Lemma \ref{Lower bound}, i.e.,
\begin{align}\label{low-g-6}
g( {{K_t}( \omega  )} )\ge5/12.
\end{align}

\begin{claim}\label{b-v-nb-r}
Any blue vertex of $K_t(\omega) $ can not incident to a blue edge.
\end{claim}
{\em Proof.} By Lemma \ref{not generalized}, a blue generalized ${K_3}$ in $K_t(\omega) $ will yield a blue $K_3$ in ${B_n(t)}$.\hfill$\Box$

\smallskip
Suppose first that there exists a {\bf blue} vertex $x$ of $K_t(\omega)$. Let $A_x\subset B_n(t)$ be the class corresponding to $x$, and let $B_m=B_n(t)-A_x$.
Clearly, $B_m$ contains no blue ${K_3}$. Moreover, $B_m$ contains no red ${K_5}$; otherwise, from Lemma \ref{not generalized}, $K_t(\omega)-x$ contains a red generalized  $K_5$, which together with $x$ would yield a red generalized $K_6$ in $K_t(\omega)$ since all edges incident to $x$ are red from Claim \ref{b-v-nb-r}. However, then $B_n(t)$ contains a red $K_6$ from Lemma \ref{not generalized}, again a contradiction. Now, it follows by Theorem \ref{01} that $e( {{B_m}} ) \le \frac{1}{2}( {1 - \frac{1}{5}}){m^2} + o( m^2)$ since $B_m$ contains neither a  blue ${K_3}$ nor a  red $K_5$, and so we have  $e( {{B_n(t)}}) \le \frac{1}{2}( {1 - \frac{1}{6}} ){n^2} + o( {{n^2}} )$ from Lemma \ref{10}. This together with the lower bound yield
$e( {{B_n(t)}}) = \frac{1}{2}\left({1 - \frac{1}{6}}\right){n^2} + o(n^2).$
Therefore, by Lemma \ref{equal}, $|A_x|=\frac{n}{6}+o(n)$ and $e(B_m)=\frac{1}{2}( {1 - \frac{1}{5}}){m^2} + o( m^2)$. It was proved in \cite[Theorem 3]{B3} that $U(m,3,3)$ (see Construction \ref{con1}) is an asymptotically extremal sequence of $RT(m,3,5,o(m))$ with weak stability property. In particular, the corresponding weighted colored complete graph $K_5(\omega)$ satisfies that all vertices are red and the red/blue edge coloring must be a {\bf pentagonlike coloring}. Thus, $\Delta(B_m, U(m,3,3))=o(m^2)$. Since $x$ is blue from the assumption, again by Claim \ref{b-v-nb-r}, we have $\Delta(B_n(6), H(n,3,3))=o(n^2).$
We will see that only when $t=6$, such a GBE-graph $B_n(6)$ can achieve the Ramsey-Tur\'{a}n density $5/12$.

Now we assume that {\bf all vertices} of ${K_t}( \omega )$ are {\bf red}. This will lead to a contradiction.

\begin{fact} \label{sp-f}
{(i)} $t \le r( {3,5} ) = 14.$

{(ii)} Each vertex $u \in {K_t}( \omega  )$ satisfies $\deg_B(u)\le4$, i.e., $u$ has at most 4 blue neighbors.

{(iii)} All edges of a red ${K_4}$ in ${K_t}( \omega )$ are half  edges.

{(iv)} There exists one half  edge of a red ${K_3}$ in ${K_t}( \omega )$.

\end{fact}
{\em Proof.}
Indeed, by Lemma \ref{not generalized}, we know that ${K_t}(\omega)$ contains neither a  blue generalized ${K_3}$ nor a  red generalized ${K_6}$. Note that all vertices of ${K_t}( \omega )$ are red, so there is no red $K_5$ in ${K_t}( \omega )$ since a red $K_5$ with some red vertex will form a red generalized $K_6$  from Definition \ref{def4}. Thus $t \le r( {3,5} ) = 14$. Similarly, (ii), (iii) and (iv) hold.\hfill$\Box$


\bigskip\noindent
{\bf Part (I)} \; $t\le6$, or $t=9,13,14$

\medskip

If $1\le t \le 5$, then the density of $B_n(t)$ is at most $\frac{1}{2}( {1 - \frac{1}{5}} )$, which contradicts (\ref{low-g-6}).
If $t = 6$, then the density can be at most $\frac{1}{2}( {1 - \frac{1}{6}} )$, which will be achieved only when all the edges are {\bf full}.
Since there is no blue ${K_3}$, we can obtain a red ${K_3}$ in ${K_6}(\omega)$, contradicting Fact \ref{sp-f} (iv).

If $t=13$, then Fact \ref{sp-f} (ii) implies that ${\deg_R}( u ) \ge 8$ for each vertex $u \in {K_t}( \omega  )$.
Thus each $u$ is red-incident to at least $3$ {\bf half}  edges; otherwise, $u$ is red-incident to at least 6 full  edges, say $\{ux_1,ux_2, \ldots ,ux_6\}$. Since there is no blue ${K_3}$, $\{x_1,\dots,x_6\}$ must contains a red $K_3$, which together with $u$ yield a red ${K_4}$ with at least 1 full edge, contradicting Fact \ref{sp-f} (iii).
Now Lemma \ref{123} implies that
$2g( {{K_t}( \omega  )} ) < 5/6,$
contradicting (\ref{low-g-6}). Similarly, $t=14$ is impossible.

Now suppose $t=9$.
Since $r( {3,4} ) = 9$ and there is no blue ${K_3}$, we have that $ {K_t}( \omega  )$ contains a red ${K_4}$ which must be {\bf half} from Fact \ref{sp-f} (iii). Suppose $X = \{x_1,x_2,x_3,x_4\}$ forms such a red ${K_4}$, and $Y = \{ {{y_1},{y_2}, \ldots ,{y_5}} \}$ consists of the remaining vertices.
\begin{claim}\label{40}
There exists one half  edge in $Y$.
\end{claim}
\noindent
{\em Proof.} On the contrary, all edges in $Y$ are full. Then $Y$ contains no red ${K_3}$ from Fact \ref{sp-f}.  Clearly $Y$ contains no blue $K_3$. It follows from Lemma \ref{jian} that $Y$ induces a \emph{pentagonlike} coloring. By symmetry, assume that $y_1y_2 \ldots y_5y_1$ is a red cycle $C_5$.

Note that each vertex has at most 4 blue neighbors from Fact \ref{sp-f} (ii) and $Y$ induces a \emph{pentagonlike} coloring, thus $\deg_B(y_i,X) \le 2$ and so $\deg_R(y_i,X ) \ge 2$ for each ${y_i} \in Y$. If there exists some ${y_j} \in Y$ such that $\deg_R(y_j,X ) = 4$, then we can easily get a red $K_4$ with a full edge, contradicting Fact \ref{sp-f} (iii).
 Thus for each ${y_i} \in Y$, $2 \le \deg_R(y_i,X ) \le 3$. Also, there are at most two vertices of $Y$ such that each of which has at least 3 red neighbors in $X$; otherwise, there would be a red $K_4$ with a full edge. Therefore,
 $10 \le {e_R}( {X,Y} )\le 12.$
Since there is no blue ${K_3}$, we have $\deg_B(x_i,Y) \le 2$ for each ${x_i} \in X$ and so ${e_B}( {X,Y} ) \le 8$. Therefore, ${e_R}( {X,Y} ) = 12$, and ${e_B}( {X,Y} ) = 8.$
 Moreover, there are exactly two vertices in $Y$ such that each of which has exactly three red neighbors in $X$. These two vertices of $Y$ could not be two consecutive vertices; otherwise we can obtain a red $K_4$ with a full edge, a contradiction.
By symmetry, we may assume that these two vertices are $y_1$ and $y_3$, and $y_1$ is red-adjacent to ${x_1},{x_2},{x_3}$ and blue-adjacent to $x_4$. Thus $y_2$ must be red-adjacent to $x_4$ since otherwise $y_1$ and $y_2$ have two red common neighbors in $X$, contradicting Fact \ref{sp-f} (iii). Similarly, $y_5$ must be red-adjacent to $x_4$.  Therefore, $x_4$ must be blue-adjacent to $y_3$ or $y_4$ since each vertex in $X$ has exactly two blue neighbors in $Y$. However, then either $\{y_1,y_3,x_4\}$ or $\{y_1,y_4,x_4\}$ will form a blue $K_3$, which  is impossible. \hfill$\Box$

\smallskip
Now we  obtain $2g( {{K_t}( \omega  )} )<5/6$ by applying Lemma \ref{1} with $q=s=1$, and $t=9$.
The proof for $t=9$ and hence {\bf Part (I)} is complete.


\begin{fact}\label{two-ind-eg}
If $u$ is red full-adjacent to $X$ with $|X|=4$, then $X$ has two independent red half edges.
\end{fact}
\noindent{\em Proof.} Clearly, $X$ contains no blue ${K_3}$. Moreover, there is no red ${K_3}$ in $X$ by Fact \ref{sp-f} (iii). Therefore, $X$ contains two independent red edges, which must be half from Fact \ref{sp-f} (iv).\hfill$\Box$

\begin{fact}\label{pent}
For a vertex $u$, if $u$ is red full-adjacent to $X$ with $|X|=5$, then the following hold.

(i) $X$ must induce a pentagonlike coloring, in which the red cycle $C_5$ must be half.

(ii) There is no other vertex which is red-adjacent to $u$.
\end{fact}
\noindent{\em Proof.} (i) Note that $X$ contains neither a  blue ${K_3}$ nor a red ${K_3}$, so the assertion follows by Lemma \ref{jian}. Moreover, all edges of the red $C_5$ in $X$ must be half from Fact \ref{sp-f} (iv). The assertion of (ii) is clear from Fact \ref{sp-f} (iii).
\hfill$\Box$

\medskip\noindent
{\bf Part (II) \; $t=7,8$}

\medskip
We first prove that for $t=7$.

\begin{claim}\label{60}
For each $u \in {K_t}( \omega  )$, $u$ is incident to at least $1$ half  edge.
\end{claim}
\noindent
{\em Proof.}~On the contrary, there exists some vertex $u$ such that all edges incident to $u$ are full.
From Fact \ref{sp-f} (ii), we have $\deg_B(u)\le 4$.
If ${\deg_B}( u ) = 4$, then the neighborhood ${N_B}( u )$ forms a red ${K_4}$. Thus all edges in ${N_B}( u )$ are half  edges from Fact \ref{sp-f} (iii), and so there are two independent half edges. We apply Lemma \ref{1} with  $s=2$, and $t=7$ to obtain
$2g( {{K_t}( \omega  )} ) \le 1 - \frac{{30}}{{30t - 40}} < 5/6.$ Suppose $\deg_B( u) = 3$.  Let
$X = N_R( u ) =\{ {x_1,x_2,x_3} \}$, and $Y= N_B( u)= \{ y_1,y_2,y_3\}.$
 Then $Y$ forms a red ${K_3}$ in which there is a  {\bf half} edge from Fact \ref{sp-f} (iv).
Moreover, from the assumption that all edges incident to $u$ are full and there is no blue $K_3$, so there must have one {red} {\bf half}  edge in $X$ from Fact \ref{sp-f} (iii). So we assume ${\deg_B}( u ) \le 2$, then ${\deg_R}( u ) \ge 4$.
By Fact \ref{two-ind-eg}, there exist two independent half edges in $N_R(u)$. For these two cases, we again obtain $2g( {{K_t}( \omega  )} ) < 5/6$.\hfill$\Box$

\smallskip
From the above claim,  we  obtain
 $
 2g( {{K_t}( \omega  )} ) \le 1 - \frac{1}{t} - \frac{1}{{2t}} < 5/6
 $ by applying Lemma \ref{123} with $\lambda= 1$ and $t=7$.
The proof for $t=7$ is complete.

We next prove that for $t=8$.
\begin{claim}\label{50}
For each $u \in {K_t}( \omega  )$, $u$ is incident to at least $1$ half  edge.
\end{claim}
\noindent{\em Proof.} ~On the contrary, there exists some vertex $u$ such that all edges incident to $u$ are full.
From Fact \ref{sp-f} (ii), we have ${\deg_B}( u ) \le 4$.
Similar as above, ${\deg_B}( u ) = 4$ is impossible.

Suppose $\deg_B( u) = 3$.  Let
$X = N_R( u ) =\{ {x_1,x_2,x_3,x_4} \}$, and $Y= N_B( u)= \{ y_1,y_2,y_3\}.$
Then $Y$ forms a red ${K_3}$ in which there is one half  edge from Fact \ref{sp-f} (iv), say $y_1y_2$.
By Fact \ref{two-ind-eg}, we may assume that $x_1x_2$ and $x_3x_4$ are two independent {\bf red} half  edges. Then Lemma \ref{1} implies that $2g( {{K_t}( \omega  )} ) \le 1 - \frac{{30}}{{30t - 20\cdot3}} = 5/6.$
 We will show the equality does not hold.
 Otherwise, all edges but $x_1x_2$, $x_3x_4$ and $y_1y_2$ in ${{K_t}( \omega  )}$ are {\bf full}. Clearly, $x_1x_3 $, $x_1x_4$, $x_2x_3$, and $x_2x_4$ must be blue from Fact \ref{sp-f} (iv). To avoid a blue ${K_3}$, so ${N_B}( {{y_i},X} ) \subseteq \{ {{x_1},{x_2}} \}$ or ${N_B}( {{y_i},X} ) \subseteq{\{ {{x_3},{x_4}} \}}$ for all ${y_i} \in Y$. By Pigeonhole Principle, there are two vertices in $Y$ such that they are red-adjacent to $\{x_3,x_4\}$ or $\{x_1,x_2\}$, contradicting Fact \ref{sp-f} (iii).

Now we suppose ${\deg_B}( u ) \le 2$, then ${\deg_R}( u ) \ge 5$.
Let $X = \{ {{x_1},\ldots ,{x_5}} \} \subseteq {N_R}( u ).$
Since all edges between $u$ and $X$ are red full  edges,  from Fact \ref{pent} (i), there exists a red cycle $C_5$ in which all edges are half.
Thus,
$2g( {{K_t}( \omega  )} ) \le 1 - \frac{{30}}{{30t - 75}} < 5/6$ from Lemma \ref{1}.\hfill$\Box$

\smallskip
From Claim \ref{50}, we  obtain
$
2g( {{K_t}( \omega  )} ) \le 1 - \frac{1}{t} - \frac{\lambda}{{2t}} < 5/6
$ by Lemma \ref{123}, completing the proof for $t=8$.

\medskip\noindent
{\bf Part (III) \; $t=10,11,12$}

\medskip

If $t=12$, then each $u \in {K_t}( \omega  )$ is incident to at least three half  edges.
Otherwise, there exists a vertex $u$ which is incident to at most two half  edges.
From Fact \ref{sp-f} (ii), we have $\deg_B(u) \le 4$, and so $\deg_R(u) \ge 7$.
Then $u$ is red full adjacent to at least $5$ vertices in $N_R(u)$, contradicting Fact \ref{pent} (ii).
Therefore, Lemma \ref{123} implies
$
2g( {{K_t}( \omega  )} ) \le 1 - \frac{1}{t} - \frac{3}{{2t}} < 5/6.
$
If $t=11$, then similarly as above,
each $u \in {K_t}( \omega  )$ is incident to at least two half  edges.
Thus, Lemma \ref{123} again implies
$2g( {{K_t}( \omega  )} ) \le 1 - \frac{1}{t} - \frac{2}{{2t}} < 5/6.$
Finally, we  prove that for $t=10$.

\begin{claim}\label{90}
For each $u \in {K_t}( \omega  )$, $u$ is incident to at least two half  edges.
\end{claim}
\noindent
{\em Proof.}~On the contrary, there is some vertex $u$ which is incident to at most one half  edge.
Together with Fact \ref{pent} (ii), $\deg_R(u) \le 5$. Moreover, from Fact \ref{sp-f} (ii), we have $\deg_B(u) \le 4$, and so $\deg_R(u) \ge 5$.
Therefore, $\deg_R(u)=5$. Let
$X = N_R( u ) =\{ {x_1,\dots,x_5} \}$, and $Y= N_B( u)= \{ y_1,\dots,y_4\}.$
Then $Y$ forms a red ${K_4}$ and all edges in $Y$ are half from Fact \ref{sp-f} (iii).

Suppose first that all edges between $u$ and $X$ are full.
Then by Fact \ref{pent} (i), there exists a red cycle $C_5$ in $X$ in which all edges are half.
Therefore, we can apply Lemma \ref{1} with $p=q=1$, and $t=10$ to obtain that
$2g( {{K_t}( \omega  )} ) \le 1 - \frac{{30}}{{30t - 75-72}} < 5/6$.
Now suppose that there is exactly one half edge between $u$ and $X$, say $ u{x_5}$.
From Fact \ref{two-ind-eg}, there exist two independent half  edges in $\{ x_1,\dots,x_4\}$. Since $Y$ forms a red $K_4$ in which all edges are half, we can apply Lemma \ref{1} with $q=1$, $s=3$ and $t=10$ to obtain
$2g( {{K_t}( \omega  )} ) \le 1 - \frac{{30}}{{30t - 72-20\cdot3}} < 5/6.$
\hfill$\Box$

\smallskip
From Claim \ref{90}, we obtain $2g( {{K_t}( \omega  )} ) \le 1 - \frac{1}{t} - \frac{2}{{2t}} < 5/6$ by Lemma \ref{123}.
The proof for $t=10$ is complete.

\medskip
From Parts (I)--(III), $\rho(3,6)=5/12$. Moreover, only when $t=6$, such a {\em GBE}-graph $B_n(6)$ can achieve the Ramsey-Tur\'{a}n density $5/12$ and $\Delta(B_n(6),H(n,3,3))=o(n^2)$.
 \hfill$\Box$

\subsection{Weak stability property for $RT(n,3,6,o(n))$}


Let $\{G_n\}$ be any sequence of asymptotically extremal graphs of the Ramsey-Tur\'{a}n problem $RT(n,3,6,o(n))$, i.e., $G_n$ is $(K_3,K_6)$-free, $\alpha ( G_n )= o( n )$, and
$e( G_n ) = RT( {n,3,6,o( n )} ) + o( n^2 )$. In order to show the weak stability property for $RT(n,3,6,o(n))$, it suffices to show that $\Delta ( G_n,H(n,3,3) ) = o( {{n^2}} )$.

For sufficiently large $n$, we take $\delta,m,\varepsilon$ and $\mu$ such that
$0<\frac{1}{n}\ll\delta\ll\frac{\mu^4}{4m^2}$, and $\frac1m\ll\varepsilon\ll\mu^2\ll 1.$
For convenience, we will also use colors {\bf1} and {\bf2} instead of {\bf blue} and {\bf red}, respectively.
For a vertex $u$ or an edge $uv$, $\chi(u)$ or $\chi(uv)$ denotes the color that $u$ or $uv$ has received.

Using Lemma \ref{b-indent} we obtain a partition $V(G_n)=\sqcup_{k=1}^2  U_k$ so that for $k\in[2]$ every $Y\subseteq U_k$ with $|Y|>\sqrt{\delta} n$ contains an edge of color $k$. Applying Lemma \ref{reg-le} to $G_n$ with $\varepsilon>0$ and $M=\frac{1}{\mu}$ we obtain an integer $M'$ and an $\varepsilon$-regular equitable partition $\sqcup_{i=0}^m V_i$ of $V(G_n)$ with $M\le m\le M'$, which refines the partition $\sqcup_{k=1}^2  U_k$. We may assume that $|V_0|=0$ and $|V_i|=\frac{n}{m}$ for $i\in[m]$ since it does not affect the result.



Now, we define a weighted and colored graph $H(\omega_0)$ on $\{v_1,v_2,\ldots,v_m\}$ (with a weight function $\omega_0$ and a distribution ${\bf{u_0}}=(u_1,u_2,\dots,u_m)$, and allowing also multiple edges), and

(i) for $k\in[2]$, we color a vertex $v_i$ by $\chi(v_i)=k$ iff $V_i\subseteq U_k$;

(ii) for $i\in[m]$, we assign the weight $u_i=\frac{|V_i|}{|V(G_n)|}=\frac{1}{m}$ to the vertex $v_i$;

(iii) for distinct $i,j\in[m]$,
we define the weight of the edge $(v_i, v_j)$ in color $k\in[2]$ as
\begin{equation}\label{weight}
\omega_{k}(v_i, v_j)=
\begin{cases}
0&\mbox {if~$d_k(V_i,V_j)<2\mu$~or~$(V_i,V_j)$~is~not~$\varepsilon$-regular,}\\
\frac{1}{2}&\mbox {if~$2\mu \le d_k(V_i,V_j)<\frac{1}{2}+2\mu$ and~$(V_i,V_j)$~is~$\varepsilon$-regular,}\\
1&\mbox {if~$\frac{1}{2}+2\mu\le d_k(V_i,V_j)$ and~$(V_i,V_j)$~is~$\varepsilon$-regular.}
\end{cases}
\end{equation}

If the ``total density''
$\sum_{k  = 1}^2 {\omega_k(v_i,v_j)}  > 1$, say $\omega_1(v_i,v_j)=\frac{1}{2}$ and $\omega_2(v_i,v_j)=1$, then we will assign $ \omega_2(v_i,v_j)=\frac{1}{2}$ to obtain $\sum_{k  = 1}^2 {\omega_k(v_i,v_j)}  = 1$. This will not affect the total number of edges of $G_n$. Let $e(H(\omega_0))=\sum_{i<j} {\sum_{k = 1}^2 {{\omega_{k }}(v_i,v_j)}}$, which is the number of weighted edges of $H(\omega_0)$.
Since $e(V_i,V_j)\le |V_i||V_j|$, we have that
 \begin{align}\label{w-eq}
 e(G_n) \le \sum\limits_{i<j} {\sum\limits_{k = 1}^2 {{\omega_{{k }}}(v_i,v_j){u_i}} {u_j}} {n^2} + 2\mu {n^2} + 4{\varepsilon}{n^2} + \frac{n^2}{m}
 \le e(H(\omega_0))\left(\frac{n}{m}\right)^2+4\mu {n^2},
 \end{align}
where the term $2\mu n^2$ comes from pairs $(V_i,V_j)$ with $d_{k}(V_i,V_j)<2\mu$ or $\frac12<d_{k}(V_i,V_j)<\frac{1}{2}+2\mu$, and $4\varepsilon n^2$ estimates the number of edges corresponding to irregular pairs, and the term $\frac{n^2}{m}<\mu n^2$ represents the internal edges lie in $V_i$.  From (\ref{w-eq}), we have $e(H(\omega_0))\ge \frac{5}{12}m^2-5\mu m^2.$
\begin{claim}\label{no K_3 nor K_6}
$H(\omega_0)$ contains neither a  blue generalized $K_3$ nor a  red generalized $K_6$.
\end{claim}
{\em Proof.}~ On the contrary, we first suppose that $H(\omega_0)$ contains a red generalized $K_6$ with vertex set $X$. Let $Y\subseteq X$ be the ``marked" subset in which all edges are full ($Y=\emptyset$ is possible). Then $Y\subseteq U_2$ from (i). For convenience, let $|X|=a$ and $|Y|=b$, then $a+b=6$.
Let $V_1,\dots, V_a$ be the parts correspond to the vertices of $X:=\{v_1,\dots,v_a\}$.
We have four cases. By a standard argument using Lemma \ref{M-s-d} and Lemma \ref{Sli-L} repeatedly we can obtain a red $K_6$ in $G_n$ for each case, and we only verify the case when $a=b=3$.

For this case, we may assume that $X=Y=\{v_1, v_2,v_3\}$, then for each $i\neq j\in[3]$, $(V_i,V_j)$ is a $(\varepsilon,\frac{1}{2}+2\mu)$-regular pair in red. From Lemma \ref{M-s-d}, there exists a subset $V_3'\subseteq V_3$ with $|V_3'|\ge(1-\sqrt{\varepsilon})|V_3|$ such that each vertex in $V_3'$ is red adjacent to at least $(\frac{1}{2}+\mu)|V_i|$ vertices in $V_i$ for $i\in[2]$. Since $|V_3'|\ge(1-\sqrt{\varepsilon})|V_3|\gg\sqrt{\delta}n$, there is a red edge, say $w_1w_2$, in $V_3'$.
 For $i\in[2]$, let $V_i'$ denote the red common neighborhood of $w_1$ and $w_2$ in $V_i$. Then $|V_i'|\ge2\mu|V_i|$ for $i\in[2]$ and $(V_1',V_2')$ is a $(\sqrt{\varepsilon},\frac{1}{2}+\mu)$-regular pair from Lemma \ref{Sli-L}.
 We now apply Lemma \ref{M-s-d} to $(V_1',V_2')$ to obtain a subset $V_1''\subseteq V_1'$ with $|V_1''|\ge (1-\sqrt{\varepsilon})|V_1'|$ such that each vertex in $V_1''$ is red adjacent to at least $(\frac{1}{2}+\mu-\sqrt{\varepsilon})|V_2'|\ge(\frac{1}{2}+\frac{\mu}{2})|V_2'|$ vertices in $V_2'$. Since $|V_1''|\ge(1-\sqrt{\varepsilon})(2\mu-\sqrt{\varepsilon})|V_1|\gg\sqrt{\delta}n$, there exists a red edge in $V_1''$, say $w_3w_4$. Note that the number of the red common neighbors of $w_3$ and $w_4$ in $V_2'$ is at least $\mu|V_2'|\gg\sqrt{\delta}n$, thus there exists a red edge in $V_2'$, denoted by $w_5w_6$, and so $\{w_1,\ldots,w_6\}$ forms a red $K_6$ in $G_n$, a contradiction.

Thus, $H(\omega_0)$ contains no red generalized $K_6$. Similarly, we can show that $H(\omega_0)$ has no blue generalized $K_3$. This completes the proof of the claim.\hfill$\Box$

\begin{definition}
Given a $2$-edge-coloring $\varphi$ of $G$, i.e., $\varphi:E(G)\rightarrow[2]$. We will write $\varphi(A,B)=k$ if $\varphi(e)=k$ for all $e\in E(G[A,B])$ and write $\varphi(v,B)$ instead of $\varphi(\{v\},B)$. If $\varphi$ is also defined on $V(G)$, we write $\varphi(A)=k$ if $\varphi(v)=k$ for all $v\in A$.
\end{definition}

We define a {\bf weighted} Tur\'{a}n graph, denoted $F_{m,6}$, as follows.
\begin{definition}\label{Fm6}
Let $T_{m,6}$ be the Tur\'{a}n graph with $6$ parts $W_1,\ldots, W_6$, where each vertex is of weight $\frac{1}{m}$, and each edge is of weight $1$.
Such a weighted Tur\'{a}n graph is denoted by $F_{m,6}$.
\end{definition}

\begin{claim}\label{max}
We have $\frac{5}{12}m^2-5\mu m^2\le e(H(\omega_0))\le \frac{5}{12}m^2.$
Furthermore, $e(H(\omega_0))= \frac{5}{12}m^2$ iff $H(\omega_0)$ isomorphic to $F_{m,6}$.
\end{claim}
{\em Proof.}~Note that $e(H(\omega_0))\ge \frac{5}{12}m^2-5\mu m^2$. In order to obtain the extremal structure of $H(\omega_0)$, we will apply the following symmetrization method to the graph $H(\omega_0)$. Suppose that $\omega_k(v_i,v_j)=0$ for each $k\in[2]$. If the weighted degrees of vertices $v_i$ and $v_j$ in $H(\omega_0)$ satisfy $d(v_i)\ge d(v_j)$, then we replace $v_j$ with a copy of $v_i$, i.e., change $\chi(v_j)$ to $\chi(v_i)$ and change $\omega_k(v_j,v_\ell)$ to $\omega_k(v_i,v_\ell)$ for $k\in[2]$ and for all $\ell\neq i,j$. Then the number of edges changes by $d(v_i)- d(v_j)\ge 0$. 
Note that this operation will not increase the maximum size of a generalized complete subgraph in any colors since $\omega_k(v_i,v_j)=0$ for each $k\in[2]$. We will repeat the operation until there is no pair $(v_i,v_j)$ with $\omega_k(v_i,v_j)=0$ for each $k\in[2]$. We now have an equivalence partition $\sqcup_{i=1}^t A_i$ for $\{v_1,\dots,v_m\}$ as follows: two vertices $v_i$ and $v_j$ are in the same class iff $\chi(v_i)=\chi(v_j)$ and $\omega_k(v_i,v_j)=0$ for each $k\in[2]$. Denote $A_i=\langle a_i\rangle$ and $a_i$ is the representative of class $A_i$. Therefore, for fixed $1\le i,j\le t$, all the edges between $A_i$ and $A_j$ have equal weights of color $k$ for $k\in[2]$, which we denote by $\omega_k(A_i,A_j)$, that is, for all vertices $u\in A_i$ and $v\in A_j$, we have $\omega_k(u,v)=\omega_k(a_i,a_j)$. Moreover, for all vertices $x\in A_i$, we have $\chi(x)=\chi(a_i)$ and $d(x)=d(a_i)$.

Now, we obtain a weighted colored complete graph $K_{t}(\omega)$ defined on $\{a_1,\ldots,a_t\}$ such that
for each $i\in [t]$ we provide $a_i$ with the weight $u_{a_i}=\sum_{v_j\in A_i} {u_j}=\frac{|A_i|}{m}$, and each edge has weight $1$ or $1/2$.
Moreover, if a pair $(a_i,a_j)$ satisfies $\omega_1(a_i,a_j)=\omega_2(a_i,a_j)=\frac{1}{2}>0$, then we replace these two edges between $a_i$ and $a_j$ by one edge of color $3-k$ and multiplicity $1$ if $\chi (a_i)=k$ where $k\in[2]$. We may assume $\chi (a_i)=2$. Then the maximum size of a generalized complete subgraph in color 2 does not increase. Moreover, the maximum size of a generalized $K_q$ in color 1 does not increase. Indeed, suppose that after this step the sets $X$ and $Y$ with $Y\subseteq X$ define a generalized $K_q$ of color 1 and $a_i$, $a_j\in X$. Then $a_i\notin Y$ since $\chi (a_i)=2$, and so the contribution of $a_i$ to $|X|+|Y|$ is 1, not depending on the weight of $(a_i,a_j)$. Clearly, the edge-density will not increase at this step.

Note that if there exists a generalized $K_p$ in $K_{t}(\omega)$, it follows from the definition of $K_{t}(\omega)$ that $H(\omega_0)$ contains a generalized $K_p$ too. Thus, $K_{t}(\omega)$ contains neither a  blue generalized $K_3$ nor a  red generalized $K_6$ from Claim \ref{no K_3 nor K_6}.

We know that $\rho(3,6)=\frac{5}{12}$ from Subsection \ref{3.2}, and only  $B_n(6)$ can achieve the Ramsey-Tur\'{a}n density $5/12$. Therefore, the corresponding weighted colored complete graph $K_6(\omega)$ defined on $\{a_1,a_2,\dots,a_6\}$ with all edges of $K_6(\omega)$ are full and the weight of each vertex in $K_6(\omega)$ equals to $\frac{1}{6}$ is the unique graph that achieves such density, and we may assume that $a_6$ is colored blue and all other vertices are colored red, and all edges $a_6a_i$ for $i\in[5]$ are colored red, and the red/blue edge coloring of $K_5$ induced by $\{a_1,\ldots,a_5\}$ is a pentagonlike coloring.

Therefore, $H(\omega_0)$ attains the maximum number of edges {\em only if} $|A_i|=\frac{m}{6}$ for $i\in[6]$, and $\chi(A_6)=\chi(a_6)=1$, and $\chi(A_i)=\chi(a_i)=2$ for $i\in[5]$, and $\chi(A_i,A_j)=\chi(a_i,a_j)$ and $\omega(x,y)=0$, for all $x,y\in A_i$ and $1\le i<j\le6$. Consequently,
$
e(H(\omega_0))\le\rho(3,6)m^2=\frac{5}{12}m^2,
$
and if $H(\omega_0)$ attains the maximum number of edges, then $H(\omega_0)$ isomorphic to $F_{m,6}$.
\hfill$\Box$

\smallskip

From the above claim, if $G_n$ is an asymptotically extremal graph of the Ramsey-Tur\'{a}n problem $RT(n,3,6,o(n))$, then $\Delta(H(w_0),F_{m,6})=o(m^2).$
Recall the definition of $H(n,3,3)$ from Lemma \ref{Lower bound}, $F_{m,6}$ can be referred to as a reduced graph of $H(n,3,3)$ by deleting the $o(n^2)$ internal edges of the six parts.
Therefore, we conclude that $\Delta ( G_n,H(n,3,3) ) = o( {{n^2}} ),$ so the weak stability holds.
The proof of Theorem \ref{zhu} is complete. \hfill$\Box$

\section{Proof of Theorem \ref{zhu-2}}\label{k3-k7}


\subsection{$\rho(3,7)=7/16$}\label{v-3-7}

From Construction \ref{con1} and $r(3,4)=9$, we know that
$ \rho(3,7)\ge \frac{1}{2}(1 - \frac{1}{r(3,4) - 1})=7/{16}
$
as desired. In the following, we will show $\rho(3,7)\le7/16$.
\medskip

Given the cyclic group $\mathbb{Z}_{3k-1}$ and a sum-free subset $S_k=\{k,k+1,\ldots,2k-1\}$, the Cayley graph $\text{Cay}_k=\text{Cay}_k(\mathbb{Z}_{3k-1},S_k)$ is defined on $\mathbb{Z}_{3k-1}$, in which for distinct $i,j\in\mathbb{Z}_{3k-1}$,
\[
ij\in E(\text{Cay}_k)\Leftrightarrow |i-j|\in S_k.
\]
Consider a complete graph $K_{3k-1}$ with $V(K_{3k-1})=\{v_0,v_1,\ldots,v_{3k-2}\}$. Call a red/blue edge coloring of $K_{3k-1}$ a \emph{nice coloring} if
 $v_iv_j \text{ is in blue}\Leftrightarrow ij\in E(\text{Cay}_k).$
When $k=2$, the nice coloring is the {\em pentagonlike coloring} which has been introduced in Section \ref{3.2}.

For $k=3$, the nice coloring of the complete graph $K_{8}$ is the eightgonlike coloring, which will be denoted by $\Gamma_{8}$. Moreover, we define $\Gamma_{8}'$ as the colored graph obtained from $\Gamma_{8}$ by recolor the blue edge $v_0v_4$ of $\Gamma_{8}$ red, and define $\Gamma_{8}''$ as the colored graph obtained from $\Gamma_{8}$ by recolor the blue edge $\{v_0v_4,v_1v_5\}$ of $\Gamma_{8}$ red.  One can see $\Gamma_{8}, \Gamma_8', \Gamma_8''$ in Figure \ref{P-E}.

Given two complete graphs $K_t$ and $K_\ell$, recall that the Ramsey number $r(t,\ell)$ is the minimum integer $N$ such that any red/blue edge coloring of the complete graph $K_N$ contains either a red $K_t$ or a blue $K_\ell$. Therefore, there exists a red/blue edge coloring of $K_{N-1}$ such that it  contains neither a red $K_t$ nor a blue $K_\ell$, for which colored complete graph is called a \textbf{Ramsey graph} for $r(t,\ell)$.

\medskip

To obtain the upper bound of $\rho(3,7)$, we need the following lemma due to K\'{e}ry \cite{kery} which tells that the Ramsey graphs for $r(4,3)$ are not unique.
\begin{lemma}[K\'{e}ry \cite{kery}]\label{p-e-like-color}
There are three Ramsey graphs for $r(4,3)$ in the sense of automorphism, i.e.,  $\Gamma_{8}, \Gamma_8'$, and $\Gamma_8''$.
\end{lemma}


\begin{figure}[h]
\begin{center}

\tikzset{every picture/.style={line width=0.75pt}} 

\tikzset{every picture/.style={line width=0.75pt}} 

\tikzset{every picture/.style={line width=0.75pt}} 

\tikzset{every picture/.style={line width=0.75pt}} 

\tikzset{every picture/.style={line width=0.75pt}} 

\tikzset{every picture/.style={line width=0.75pt}} 

\tikzset{every picture/.style={line width=0.75pt}} 

\tikzset{every picture/.style={line width=0.75pt}} 

\tikzset{every picture/.style={line width=0.75pt}} 

\tikzset{every picture/.style={line width=0.75pt}} 

\tikzset{every picture/.style={line width=0.75pt}} 

\tikzset{every picture/.style={line width=0.75pt}} 

\begin{tikzpicture}[x=0.75pt,y=0.75pt,yscale=-1,xscale=1]

\draw [color={rgb, 255:red, 208; green, 2; blue, 27 }  ,draw opacity=1 ]   (503.12,90.13) -- (503.12,219.35) ;
\draw [color={rgb, 255:red, 208; green, 2; blue, 27 }  ,draw opacity=1 ]   (459.77,109.05) -- (546.46,200.43) ;
\draw [color={rgb, 255:red, 208; green, 2; blue, 27 }  ,draw opacity=1 ]   (323.12,88.13) -- (323.12,217.35) ;
\draw  [color={rgb, 255:red, 208; green, 2; blue, 27 }  ,draw opacity=1 ] (204.42,153.74) -- (186.46,199.43) -- (143.12,218.35) -- (99.77,199.43) -- (81.81,153.74) -- (99.77,108.05) -- (143.12,89.13) -- (186.46,108.05) -- cycle ;
\draw [color={rgb, 255:red, 208; green, 2; blue, 27 }  ,draw opacity=1 ]   (98.99,109.5) -- (99.77,199.43) ;
\draw [color={rgb, 255:red, 208; green, 2; blue, 27 }  ,draw opacity=1 ]   (81.81,153.74) -- (143.12,218.35) ;
\draw [color={rgb, 255:red, 208; green, 2; blue, 27 }  ,draw opacity=1 ]   (143.12,218.35) -- (204.42,153.74) ;
\draw [color={rgb, 255:red, 208; green, 2; blue, 27 }  ,draw opacity=1 ]   (99.77,108.05) -- (186.46,108.05) ;
\draw [color={rgb, 255:red, 208; green, 2; blue, 27 }  ,draw opacity=1 ]   (143.12,89.13) -- (81.81,153.74) ;
\draw [color={rgb, 255:red, 208; green, 2; blue, 27 }  ,draw opacity=1 ]   (186.46,199.43) -- (137.85,199.43) -- (99.77,199.43) ;
\draw [color={rgb, 255:red, 208; green, 2; blue, 27 }  ,draw opacity=1 ]   (143.12,89.13) -- (204.42,153.74) ;
\draw [color={rgb, 255:red, 208; green, 2; blue, 27 }  ,draw opacity=1 ]   (186.46,108.05) -- (186.46,199.43) ;
\draw  [fill={rgb, 255:red, 0; green, 0; blue, 0 }  ,fill opacity=1 ] (78.82,153.74) .. controls (78.82,152) and (80.16,150.59) .. (81.81,150.59) .. controls (83.46,150.59) and (84.8,152) .. (84.8,153.74) .. controls (84.8,155.48) and (83.46,156.89) .. (81.81,156.89) .. controls (80.16,156.89) and (78.82,155.48) .. (78.82,153.74) -- cycle ;
\draw  [fill={rgb, 255:red, 0; green, 0; blue, 0 }  ,fill opacity=1 ] (183.47,108.05) .. controls (183.47,106.31) and (184.81,104.9) .. (186.46,104.9) .. controls (188.12,104.9) and (189.45,106.31) .. (189.45,108.05) .. controls (189.45,109.79) and (188.12,111.21) .. (186.46,111.21) .. controls (184.81,111.21) and (183.47,109.79) .. (183.47,108.05) -- cycle ;
\draw  [fill={rgb, 255:red, 0; green, 0; blue, 0 }  ,fill opacity=1 ] (96.78,199.43) .. controls (96.78,197.69) and (98.12,196.28) .. (99.77,196.28) .. controls (101.42,196.28) and (102.76,197.69) .. (102.76,199.43) .. controls (102.76,201.17) and (101.42,202.58) .. (99.77,202.58) .. controls (98.12,202.58) and (96.78,201.17) .. (96.78,199.43) -- cycle ;
\draw  [fill={rgb, 255:red, 0; green, 0; blue, 0 }  ,fill opacity=1 ] (201.43,153.74) .. controls (201.43,152) and (202.77,150.59) .. (204.42,150.59) .. controls (206.07,150.59) and (207.41,152) .. (207.41,153.74) .. controls (207.41,155.48) and (206.07,156.89) .. (204.42,156.89) .. controls (202.77,156.89) and (201.43,155.48) .. (201.43,153.74) -- cycle ;
\draw  [fill={rgb, 255:red, 0; green, 0; blue, 0 }  ,fill opacity=1 ] (140.12,218.35) .. controls (140.12,216.61) and (141.46,215.2) .. (143.12,215.2) .. controls (144.77,215.2) and (146.11,216.61) .. (146.11,218.35) .. controls (146.11,220.09) and (144.77,221.51) .. (143.12,221.51) .. controls (141.46,221.51) and (140.12,220.09) .. (140.12,218.35) -- cycle ;
\draw  [fill={rgb, 255:red, 0; green, 0; blue, 0 }  ,fill opacity=1 ] (140.12,89.13) .. controls (140.12,87.39) and (141.46,85.98) .. (143.12,85.98) .. controls (144.77,85.98) and (146.11,87.39) .. (146.11,89.13) .. controls (146.11,90.87) and (144.77,92.28) .. (143.12,92.28) .. controls (141.46,92.28) and (140.12,90.87) .. (140.12,89.13) -- cycle ;
\draw  [fill={rgb, 255:red, 0; green, 0; blue, 0 }  ,fill opacity=1 ] (183.47,199.43) .. controls (183.47,197.69) and (184.81,196.28) .. (186.46,196.28) .. controls (188.12,196.28) and (189.45,197.69) .. (189.45,199.43) .. controls (189.45,201.17) and (188.12,202.58) .. (186.46,202.58) .. controls (184.81,202.58) and (183.47,201.17) .. (183.47,199.43) -- cycle ;
\draw  [fill={rgb, 255:red, 0; green, 0; blue, 0 }  ,fill opacity=1 ] (96.78,108.05) .. controls (96.78,106.31) and (98.12,104.9) .. (99.77,104.9) .. controls (101.42,104.9) and (102.76,106.31) .. (102.76,108.05) .. controls (102.76,109.79) and (101.42,111.21) .. (99.77,111.21) .. controls (98.12,111.21) and (96.78,109.79) .. (96.78,108.05) -- cycle ;
\draw  [color={rgb, 255:red, 208; green, 2; blue, 27 }  ,draw opacity=1 ] (384.42,152.74) -- (366.46,198.43) -- (323.12,217.35) -- (279.77,198.43) -- (261.81,152.74) -- (279.77,107.05) -- (323.12,88.13) -- (366.46,107.05) -- cycle ;
\draw [color={rgb, 255:red, 208; green, 2; blue, 27 }  ,draw opacity=1 ]   (278.99,108.5) -- (279.77,198.43) ;
\draw [color={rgb, 255:red, 208; green, 2; blue, 27 }  ,draw opacity=1 ]   (261.81,152.74) -- (323.12,217.35) ;
\draw [color={rgb, 255:red, 208; green, 2; blue, 27 }  ,draw opacity=1 ]   (323.12,217.35) -- (384.42,152.74) ;
\draw [color={rgb, 255:red, 208; green, 2; blue, 27 }  ,draw opacity=1 ]   (279.77,107.05) -- (366.46,107.05) ;
\draw [color={rgb, 255:red, 208; green, 2; blue, 27 }  ,draw opacity=1 ]   (323.12,88.13) -- (261.81,152.74) ;
\draw [color={rgb, 255:red, 208; green, 2; blue, 27 }  ,draw opacity=1 ]   (366.46,198.43) -- (317.85,198.43) -- (279.77,198.43) ;
\draw [color={rgb, 255:red, 208; green, 2; blue, 27 }  ,draw opacity=1 ]   (323.12,88.13) -- (384.42,152.74) ;
\draw [color={rgb, 255:red, 208; green, 2; blue, 27 }  ,draw opacity=1 ]   (366.46,107.05) -- (366.46,198.43) ;
\draw  [fill={rgb, 255:red, 0; green, 0; blue, 0 }  ,fill opacity=1 ] (258.82,152.74) .. controls (258.82,151) and (260.16,149.59) .. (261.81,149.59) .. controls (263.46,149.59) and (264.8,151) .. (264.8,152.74) .. controls (264.8,154.48) and (263.46,155.89) .. (261.81,155.89) .. controls (260.16,155.89) and (258.82,154.48) .. (258.82,152.74) -- cycle ;
\draw  [fill={rgb, 255:red, 0; green, 0; blue, 0 }  ,fill opacity=1 ] (363.47,107.05) .. controls (363.47,105.31) and (364.81,103.9) .. (366.46,103.9) .. controls (368.12,103.9) and (369.45,105.31) .. (369.45,107.05) .. controls (369.45,108.79) and (368.12,110.21) .. (366.46,110.21) .. controls (364.81,110.21) and (363.47,108.79) .. (363.47,107.05) -- cycle ;
\draw  [fill={rgb, 255:red, 0; green, 0; blue, 0 }  ,fill opacity=1 ] (276.78,198.43) .. controls (276.78,196.69) and (278.12,195.28) .. (279.77,195.28) .. controls (281.42,195.28) and (282.76,196.69) .. (282.76,198.43) .. controls (282.76,200.17) and (281.42,201.58) .. (279.77,201.58) .. controls (278.12,201.58) and (276.78,200.17) .. (276.78,198.43) -- cycle ;
\draw  [fill={rgb, 255:red, 0; green, 0; blue, 0 }  ,fill opacity=1 ] (381.43,152.74) .. controls (381.43,151) and (382.77,149.59) .. (384.42,149.59) .. controls (386.07,149.59) and (387.41,151) .. (387.41,152.74) .. controls (387.41,154.48) and (386.07,155.89) .. (384.42,155.89) .. controls (382.77,155.89) and (381.43,154.48) .. (381.43,152.74) -- cycle ;
\draw  [fill={rgb, 255:red, 0; green, 0; blue, 0 }  ,fill opacity=1 ] (320.12,217.35) .. controls (320.12,215.61) and (321.46,214.2) .. (323.12,214.2) .. controls (324.77,214.2) and (326.11,215.61) .. (326.11,217.35) .. controls (326.11,219.09) and (324.77,220.51) .. (323.12,220.51) .. controls (321.46,220.51) and (320.12,219.09) .. (320.12,217.35) -- cycle ;
\draw  [fill={rgb, 255:red, 0; green, 0; blue, 0 }  ,fill opacity=1 ] (320.12,88.13) .. controls (320.12,86.39) and (321.46,84.98) .. (323.12,84.98) .. controls (324.77,84.98) and (326.11,86.39) .. (326.11,88.13) .. controls (326.11,89.87) and (324.77,91.28) .. (323.12,91.28) .. controls (321.46,91.28) and (320.12,89.87) .. (320.12,88.13) -- cycle ;
\draw  [fill={rgb, 255:red, 0; green, 0; blue, 0 }  ,fill opacity=1 ] (363.47,198.43) .. controls (363.47,196.69) and (364.81,195.28) .. (366.46,195.28) .. controls (368.12,195.28) and (369.45,196.69) .. (369.45,198.43) .. controls (369.45,200.17) and (368.12,201.58) .. (366.46,201.58) .. controls (364.81,201.58) and (363.47,200.17) .. (363.47,198.43) -- cycle ;
\draw  [fill={rgb, 255:red, 0; green, 0; blue, 0 }  ,fill opacity=1 ] (276.78,107.05) .. controls (276.78,105.31) and (278.12,103.9) .. (279.77,103.9) .. controls (281.42,103.9) and (282.76,105.31) .. (282.76,107.05) .. controls (282.76,108.79) and (281.42,110.21) .. (279.77,110.21) .. controls (278.12,110.21) and (276.78,108.79) .. (276.78,107.05) -- cycle ;
\draw  [color={rgb, 255:red, 208; green, 2; blue, 27 }  ,draw opacity=1 ] (564.42,154.74) -- (546.46,200.43) -- (503.12,219.35) -- (459.77,200.43) -- (441.81,154.74) -- (459.77,109.05) -- (503.12,90.13) -- (546.46,109.05) -- cycle ;
\draw [color={rgb, 255:red, 208; green, 2; blue, 27 }  ,draw opacity=1 ]   (458.99,110.5) -- (459.77,200.43) ;
\draw [color={rgb, 255:red, 208; green, 2; blue, 27 }  ,draw opacity=1 ]   (441.81,154.74) -- (503.12,219.35) ;
\draw [color={rgb, 255:red, 208; green, 2; blue, 27 }  ,draw opacity=1 ]   (503.12,219.35) -- (564.42,154.74) ;
\draw [color={rgb, 255:red, 208; green, 2; blue, 27 }  ,draw opacity=1 ]   (459.77,109.05) -- (546.46,109.05) ;
\draw [color={rgb, 255:red, 208; green, 2; blue, 27 }  ,draw opacity=1 ]   (503.12,90.13) -- (441.81,154.74) ;
\draw [color={rgb, 255:red, 208; green, 2; blue, 27 }  ,draw opacity=1 ]   (546.46,200.43) -- (497.85,200.43) -- (459.77,200.43) ;
\draw [color={rgb, 255:red, 208; green, 2; blue, 27 }  ,draw opacity=1 ]   (503.12,90.13) -- (564.42,154.74) ;
\draw [color={rgb, 255:red, 208; green, 2; blue, 27 }  ,draw opacity=1 ]   (546.46,109.05) -- (546.46,200.43) ;
\draw  [fill={rgb, 255:red, 0; green, 0; blue, 0 }  ,fill opacity=1 ] (438.82,154.74) .. controls (438.82,153) and (440.16,151.59) .. (441.81,151.59) .. controls (443.46,151.59) and (444.8,153) .. (444.8,154.74) .. controls (444.8,156.48) and (443.46,157.89) .. (441.81,157.89) .. controls (440.16,157.89) and (438.82,156.48) .. (438.82,154.74) -- cycle ;
\draw  [fill={rgb, 255:red, 0; green, 0; blue, 0 }  ,fill opacity=1 ] (543.47,109.05) .. controls (543.47,107.31) and (544.81,105.9) .. (546.46,105.9) .. controls (548.12,105.9) and (549.45,107.31) .. (549.45,109.05) .. controls (549.45,110.79) and (548.12,112.21) .. (546.46,112.21) .. controls (544.81,112.21) and (543.47,110.79) .. (543.47,109.05) -- cycle ;
\draw  [fill={rgb, 255:red, 0; green, 0; blue, 0 }  ,fill opacity=1 ] (456.78,200.43) .. controls (456.78,198.69) and (458.12,197.28) .. (459.77,197.28) .. controls (461.42,197.28) and (462.76,198.69) .. (462.76,200.43) .. controls (462.76,202.17) and (461.42,203.58) .. (459.77,203.58) .. controls (458.12,203.58) and (456.78,202.17) .. (456.78,200.43) -- cycle ;
\draw  [fill={rgb, 255:red, 0; green, 0; blue, 0 }  ,fill opacity=1 ] (561.43,154.74) .. controls (561.43,153) and (562.77,151.59) .. (564.42,151.59) .. controls (566.07,151.59) and (567.41,153) .. (567.41,154.74) .. controls (567.41,156.48) and (566.07,157.89) .. (564.42,157.89) .. controls (562.77,157.89) and (561.43,156.48) .. (561.43,154.74) -- cycle ;
\draw  [fill={rgb, 255:red, 0; green, 0; blue, 0 }  ,fill opacity=1 ] (500.12,219.35) .. controls (500.12,217.61) and (501.46,216.2) .. (503.12,216.2) .. controls (504.77,216.2) and (506.11,217.61) .. (506.11,219.35) .. controls (506.11,221.09) and (504.77,222.51) .. (503.12,222.51) .. controls (501.46,222.51) and (500.12,221.09) .. (500.12,219.35) -- cycle ;
\draw  [fill={rgb, 255:red, 0; green, 0; blue, 0 }  ,fill opacity=1 ] (500.12,90.13) .. controls (500.12,88.39) and (501.46,86.98) .. (503.12,86.98) .. controls (504.77,86.98) and (506.11,88.39) .. (506.11,90.13) .. controls (506.11,91.87) and (504.77,93.28) .. (503.12,93.28) .. controls (501.46,93.28) and (500.12,91.87) .. (500.12,90.13) -- cycle ;
\draw  [fill={rgb, 255:red, 0; green, 0; blue, 0 }  ,fill opacity=1 ] (543.47,200.43) .. controls (543.47,198.69) and (544.81,197.28) .. (546.46,197.28) .. controls (548.12,197.28) and (549.45,198.69) .. (549.45,200.43) .. controls (549.45,202.17) and (548.12,203.58) .. (546.46,203.58) .. controls (544.81,203.58) and (543.47,202.17) .. (543.47,200.43) -- cycle ;
\draw  [fill={rgb, 255:red, 0; green, 0; blue, 0 }  ,fill opacity=1 ] (456.78,109.05) .. controls (456.78,107.31) and (458.12,105.9) .. (459.77,105.9) .. controls (461.42,105.9) and (462.76,107.31) .. (462.76,109.05) .. controls (462.76,110.79) and (461.42,112.21) .. (459.77,112.21) .. controls (458.12,112.21) and (456.78,110.79) .. (456.78,109.05) -- cycle ;

\draw (131.66,251) node [anchor=north west][inner sep=0.75pt]    {$\Gamma _{8}$};
\draw (135,71.05) node [anchor=north west][inner sep=0.75pt]    {$v_{0}$};
\draw (84,91.05) node [anchor=north west][inner sep=0.75pt]    {$v_{1}$};
\draw (61,149.05) node [anchor=north west][inner sep=0.75pt]    {$v_{2}$};
\draw (80,198.05) node [anchor=north west][inner sep=0.75pt]    {$v_{3}$};
\draw (132,223.05) node [anchor=north west][inner sep=0.75pt]    {$v_{4}$};
\draw (189.45,199.43) node [anchor=north west][inner sep=0.75pt]    {$v_{5}$};
\draw (207.41,153.74) node [anchor=north west][inner sep=0.75pt]    {$v_{6}$};
\draw (186,93.05) node [anchor=north west][inner sep=0.75pt]    {$v_{7}$};
\draw (311.66,250) node [anchor=north west][inner sep=0.75pt]    {$\Gamma '_{8}$};
\draw (315,70.05) node [anchor=north west][inner sep=0.75pt]    {$v_{0}$};
\draw (264,90.05) node [anchor=north west][inner sep=0.75pt]    {$v_{1}$};
\draw (241,148.05) node [anchor=north west][inner sep=0.75pt]    {$v_{2}$};
\draw (260,197.05) node [anchor=north west][inner sep=0.75pt]    {$v_{3}$};
\draw (312,222.05) node [anchor=north west][inner sep=0.75pt]    {$v_{4}$};
\draw (369.45,198.43) node [anchor=north west][inner sep=0.75pt]    {$v_{5}$};
\draw (387.41,152.74) node [anchor=north west][inner sep=0.75pt]    {$v_{6}$};
\draw (366,92.05) node [anchor=north west][inner sep=0.75pt]    {$v_{7}$};
\draw (491.66,252) node [anchor=north west][inner sep=0.75pt]    {$\Gamma ''_{8}$};
\draw (495,72.05) node [anchor=north west][inner sep=0.75pt]    {$v_{0}$};
\draw (444,92.05) node [anchor=north west][inner sep=0.75pt]    {$v_{1}$};
\draw (421,150.05) node [anchor=north west][inner sep=0.75pt]    {$v_{2}$};
\draw (440,199.05) node [anchor=north west][inner sep=0.75pt]    {$v_{3}$};
\draw (492,224.05) node [anchor=north west][inner sep=0.75pt]    {$v_{4}$};
\draw (549.45,200.43) node [anchor=north west][inner sep=0.75pt]    {$v_{5}$};
\draw (567.41,154.74) node [anchor=north west][inner sep=0.75pt]    {$v_{6}$};
\draw (546,94.05) node [anchor=north west][inner sep=0.75pt]    {$v_{7}$};

\end{tikzpicture}

\begin{center}
\begin{caption}
{\text{Ramsey graphs for $r(4,3)$,}
\text{where all blue edges are omitted.} }\label{P-E}
\end{caption}
\end{center}
\end{center}
\end{figure}

We also need the following lemma.
\begin{lemma}\label{1-sys}
Let, in ${K_t}( \omega )$, $\ell_1$ disjoint $C_8$'s, $\ell_2$ disjoint ${K_5}$'s, $\ell_3$ disjoint ${C_5}$'s, $\ell_4$ disjoint ${K_4}$'s, $\ell_5$ disjoint ${C_4}$'s, $\ell_6$ disjoint ${K_3}$'s and $\ell_7$ independent edges be given, where these $\ell_7$ independent edges and $C_8$'s, ${K_5}$'s, ${K_4}$'s, ${C_4}$'s and ${K_3}$'s are also pairwise vertex disjoint and all are half. Then
$2g( {{K_t}( \omega  )} ) \le 1 - \frac{{30}}{{30t -120\ell_1- 100\ell_2 - 75\ell_3-72\ell_4 - 60\ell_5-45\ell_6 - 20\ell_7}}.$
\end{lemma}
\noindent
{\em Proof.} Let ${K_t}( \omega )$ be defined on $[t]$. We may assume that all of the other edges are full. There are $\ell_8 = t -8\ell_1- 5\ell_2-5\ell_3-4\ell_4-4\ell_5-3\ell_6-2\ell_7 $ ``other vertices". This structure will be denoted by $W( {\ell_1,\ldots,\ell_8} )$. By a similar argument as Lemma \ref{1}, we have $W:=W( {\ell_1,\ldots,\ell_8} )$ is a dense structure.

Let the optimum distribution of $W$ be $z_1, \ldots ,z_1$, $z_2, \ldots ,z_2$, $z_3, \ldots ,z_3$, $z_4, \ldots ,z_4$, $z_5, \ldots ,z_5$, $z_6, \ldots ,z_6$, $z_7, \ldots ,z_7$, $z_8, \ldots ,z_8$, then $$8\ell_1z_1+5\ell_2z_2+5\ell_3z_3+4\ell_4z_4+4\ell_5z_5+3\ell_6z_6+2\ell_7z_7+\ell_8z_8=1$$ due to $\sum_{i = 1}^t {{u_i}}  = 1$. Again by a similar argument as Lemma \ref{1}, the weighted degrees of the vertices of $C_8$, $K_5$, $C_5$, $K_4$, $C_4$, $K_3$, $K_2$, and
``other vertices'' are $t(1 - 2z_1)$, $t(1-3z_2)$, $t(1 - 2z_3)$, $t(1-{5z_4}/2)$, $t(1 - 2z_5)$, $t(1 - 2z_6)$, $t(1-{3z_7}/2)$, and $t(1-z_8)$, respectively. By Lemma \ref{0}, all the weighted degrees are equal.
 Hence, if we set $z_1=15m$, then $z_2= 10m$, $z_3= 15m$, $z_4= 12m$, $z_5= 15m$, $z_6= 15m$, $z_7= 20m$, and $z_8 = 30m$. Therefore,
 \begin{align*}
   1 &= 8\ell_1\cdot 15m+5\ell_2\cdot10m + 5\ell_3\cdot15m+ 4\ell_4\cdot12m + 4\ell_5\cdot15m+ 3\ell_6\cdot15m + 2\ell_7\cdot20m
   \\& \quad \quad +   ( {t -8\ell_1- 5\ell_2-5\ell_3-4\ell_4-4\ell_5-3\ell_6-2\ell_7} )\cdot30m
    \\&= ( {30t -120\ell_1- 100\ell_2 - 75\ell_3-72\ell_4-60\ell_5-45\ell_6-20\ell_7} )m.
\end{align*}
By Lemma \ref{0},
$2g( {{K_t}( \omega  )} ) \le 1 - \frac{{30}}{{30t -120\ell_1- 100\ell_2 - 75\ell_3-72\ell_4 - 60\ell_5-45\ell_6 - 20\ell_7}}$
as desired.\hfill$\Box$

\medskip

In the following, we will prove the upper bound $\rho(3,7)\le7/16$.
From Theorem \ref{00}, let $B_n(t)=B(h, t |{n_1}, \ldots , {{n_t}} |\mu)$, which is a dense asymptotically extremal graph of $RT( {n,{3},{7},o( n )} )$, with its canonical coloring with respect to the $t$ classes $X_1,\dots,X_t$ of $B_n(t)$. Let ${K_t}( \omega  )$ be the corresponding weighted colored complete graph with distribution ${\bf{u}}=(u_1,\dots,u_t)$ \;$( {u_i} > 0$, $\sum_{i = 1}^t {{u_i}}  = 1)$. We say a vertex of ${K_t}( \omega  )$ is blue (red) if all edges in the corresponding set $X_\ell$ for some $\ell\in[t]$ are blue (red).
From the assumption, the colored $B_n(t)$ contains neither a  blue $K_3$ nor a  red $K_7$, and $\alpha(B_n(t))=o(n)$. We may assume that $e(B_n(t))\ge\frac{7}{16}{n^2} + o(n^2)$ from the lower bound that $\rho(3,7)\ge7/16$, i.e., $g( {{K_t}( \omega  )} )\ge7/16.$

\begin{claim}\label{b-v-nb-r-for-K3-K7}
Any blue vertex of $K_t(\omega) $ can not incident to a blue edge.
\end{claim}

Suppose first that there exists a {\bf blue} vertex $x$ of $K_t(\omega)$. Let $X\subset B_n(t)$ be the class corresponding to $x$, and let $B_m=B_n(t)-X$.
We claim that $B_m$ contains no red ${K_6}$. Otherwise, from Lemma \ref{not generalized}, $K_t(\omega)-x$ contains a red generalized $K_6$, which together with $x$ yield a red generalized $K_7$ in $K_t(\omega)$ from Claim \ref{b-v-nb-r-for-K3-K7}. Thus, $B_n(t)$ contains a red $K_7$ from Lemma \ref{not generalized} again, a contradiction. Since $B_m$ contains no blue ${K_3}$, Theorem \ref{zhu} implies $e( {{B_m}} ) \le \frac{1}{2}( {1 - \frac{1}{6}}){m^2} + o( m^2)$, and so $e( {{B_n(t)}}) \le \frac{1}{2}( {1 - \frac{1}{7}} ){n^2} + o( {{n^2}} )$ from Lemma \ref{10},  which is impossible.

\smallskip
In the following, we assume that {\bf all vertices} of ${K_t}( \omega )$ are {\bf red}.

\begin{fact}\label{center-fact} 
(i) $t \le r( {3,6} ) = 18.$

(ii) Each vertex $u \in {K_t}( \omega  )$ satisfies $\deg_B(u)\le5$, i.e., $u$ has at most $5$ blue neighbors.

(iii) All edges of a red ${K_5}$ are half.

(iv) For each red $K_4$, it contains no full  $K_3$. This implies that it contains either a half  $K_3$ or two independent half  edges. Furthermore, if a red $K_4$ contains two independent full  edges, then it also contains two independent half  edges.

\end{fact}
{\em Proof.}
By Lemma \ref{not generalized}, we know that ${K_t}( \omega  )$ contains neither a blue generalized ${K_3}$ nor a red generalized ${K_7}$. Note that all vertices of ${K_t}( \omega )$ are red and $N_B(u)$ forms a red clique, so (i) holds from Definition \ref{def4}. Similarly, (ii), (iii), and the first assertion of (iv) hold. Therefore, for each red $K_4$, there are at least $2$ half  edges. So the remaining assertions of (iv) hold. \hfill$\Box$


\bigskip\noindent
{\bf Part (I) \; $t\le 8$ or $t=17,18$}
\medskip

If $1\le t\le7$, then the density of $B_n(t)$ is at most $\frac{1}{2}(1-\frac{1}{7})$, which is impossible.
If $t = 8$, then the density can be at most $\frac{7}{16}$, which will be achieved only when all the edges are {\bf full}, together with Fact \ref{center-fact} (iv), we have ${K_8}( \omega  )$ contains neither a red $K_4$ nor a blue ${K_3}$. Thus, the red/blue edge coloring graph of $K_8(\omega)$ must be one of the three Ramsey graphs $\{\Gamma_{8}, \Gamma_8', \Gamma_8''\}$ for $r(4,3)$ from Lemma \ref{p-e-like-color}, which implies that $\Delta(B_n(8),U(n,4,3))=o(n^2)$. We will see only $t=8$ can achieve the Ramsey-Tur\'{a}n density $7/16$.

For $t=17$, Fact \ref{center-fact} (ii) implies that ${\deg_R}( u ) \ge 11$ for each $u \in {K_t}( \omega  )$.
Clearly, each vertex $u \in K_t(\omega)$ is red-incident to at least $3$ half edges; otherwise, some vertex $u$ is red-incident to at least $9$ full edges. Since $r(3,4)=9$, and there is no blue ${K_3}$, $N_R(u)$ must contain a red $K_4$, which together with $u$ form a red ${K_5}$ with a full edge, contradicting Fact \ref{center-fact} (iii).
Thus, from Lemma \ref{123}, $2g( {{K_t}( \omega  )} ) \le 1 - \frac{1}{t} - \frac{{3}}{{2t}} < 7/8$.
$t=18$ is similar. The proof of {\bf Part (I)} is complete.


\begin{fact}\label{T-I-E}
Let $X$ be a subset of $K_t(w)$ of size $4$, if $X$ contains no red $K_3$, then $X$ contains two independent red edges.
\end{fact}

\begin{fact}\label{two-ind-eg-for-K3-K7}
If  $u\in{K_t}(\omega) $ is red full-adjacent to $X$ with $|X|=6$, then $X$ has a red half  $K_3$.

\end{fact}
\noindent{\em Proof.} Since there is no blue $K_3$, $X$ contains a red $K_3$ which must be half from Fact \ref{center-fact} (iv).\hfill$\Box$

\begin{fact}\label{eightdagon-like-K3-K7}
Suppose that $u\in{K_t}(\omega) $ is red adjacent to $X$ with $|X|=8$.
If all but at most $1$ vertices in $X$ are full-adjacent to $u$, then the following properties hold.

(i) $X$ forms one of the Ramsey graphs $\{\Gamma_{8}, \Gamma_8', \Gamma_8''\}$, and all red edges in the full-neighborhood of $u$ are half. Also, $X\cup \{u\}$ contains two disjoint red half $K_3$'s and an independent half edge.

(ii) There is no other vertex which is red-adjacent to $u$.
\end{fact}

\noindent{\em Proof.} (i) Since all but at most $1$ vertex in $X$ are full-adjacent to $u$, $X$ contains no  red ${K_4}$ from Fact \ref{center-fact} (iii). Also, $X$ has no blue ${K_3}$. Thus, $X$ forms one of the Ramsey graphs $\{\Gamma_{8}, \Gamma_8', \Gamma_8''\}$ from Lemma \ref{p-e-like-color}. Moreover, from Fact \ref{center-fact} (iv), all red edges contained in the full-neighborhood of $u$ are half since each red edge in $X$ belongs to some red $K_3$. Therefore, $X\cup \{u\}$ contains two disjoint red half $K_3$'s and an independent half edge as desired.

(ii) Suppose $v\not\in X$ is red-adjacent to $u$, then $X\cup \{v\}$ contains a red $K_4$ since there is no blue $K_3$, which together with $u$ yield a red $K_5$ with a full  edge, contradicting Fact \ref{center-fact} (iii).
\hfill$\Box$

\medskip\noindent
{\bf Part (II) \; $t=9$ or $13\le t\le16$}
\medskip

We first prove that for $t=9$.
Since there is no blue ${K_3}$, ${K_t}( \omega  )$ contains a red $K_4$ in which there exists either a red half  $K_3$ or two independent half  edges from Fact \ref{center-fact} (iv). If the former is true, we apply Lemma \ref{1-sys} to obtain $2g( {{K_t}( \omega  )} ) \le 1 - \frac{{30}}{{30t - 45}} < 7/8$. If the later occurs, we again apply Lemma \ref{1-sys} to obtain $2g( {{K_t}( \omega  )} ) \le 1 - \frac{{30}}{{30t - 40}} < 7/8$.

For $t=15,16$, we claim that
each $u \in {K_t}( \omega  )$ is incident to at least $3$ half  edges.
Otherwise, there is some vertex $u$ which is incident to at most $2$ half edges, then
Fact \ref{center-fact} (ii) implies ${\deg_B}( u ) \le 5$ since all vertices are red, and so ${\deg_R}( u ) \ge 9$. This contradicts Fact \ref{eightdagon-like-K3-K7} (ii).
Therefore, we obtain  $2g( {{K_t}( \omega  )} ) \le 1 - \frac{1}{t} - \frac{3}{{2t}} < 7/8$ by Lemma \ref{123} with $t=15$ or $16$.

We next prove that for $t=14$.
\begin{claim}\label{t=14}
For each $u \in {K_t}( \omega  )$, $u$ is incident to at least $2$ half  edges.
\end{claim}
\noindent{\em Proof.} ~On the contrary, some vertex $u \in {K_t}( \omega  )$ is incident to at most one half  edge.
From Fact \ref{center-fact} (ii), we have ${\deg_B}( u ) \le 5$, and so ${\deg_R}( u ) \ge 8$. Thus ${\deg_R}( u ) = 8$ from Fact \ref{eightdagon-like-K3-K7} (ii).
Let $X=N_R(u)$, and $Y=N_B(u).$ From Fact \ref{center-fact} (iii), $Y$ forms a red half  ${K_5}$. From Fact \ref{eightdagon-like-K3-K7} (i), $X$ forms one of the Ramsey graphs $\{\Gamma_{8}, \Gamma_8', \Gamma_8''\}$,
and $X\cup \{u\}$ contains two disjoint red half $K_3$'s and an independent half edge. We can apply Lemma \ref{1-sys} with $\ell_2=\ell_7=1$, $\ell_6=2$, and $t=14$ to obtain
$2g( {{K_t}( \omega  )} ) \le 1 - \frac{{30}}{{30t - 100-45\cdot2-20}} < 7/8.$\hfill$\Box$

\smallskip
From the above claim, we obtain
$
 2g( {{K_t}( \omega  )} ) \le 1 - \frac{1}{t} - \frac{2}{{2t}} < 7/8
$ from Lemma \ref{123}.
The proof for $t=14$ is complete.

Finally, we prove that for $t=13$.

\begin{claim}\label{t=13}
For each $u \in {K_t}( \omega  )$, $u$ is incident to at least $2$ half  edges.
\end{claim}
\noindent{\em Proof.} ~On the contrary, there is a vertex $u$ which is incident to at most $1$ half  edge. From Fact \ref{center-fact} (ii), we have ${\deg_B}( u ) \le 5$. If ${\deg_B}( u ) \le 3$, then ${\deg_R}( u ) \ge 9$, contradicting Fact \ref{eightdagon-like-K3-K7} (ii). Thus ${\deg_B}( u )= 4$, or 5.

Suppose first ${\deg_B}( u ) = 4$. Let $X=N_R(u)=\{x_1,\ldots,x_8\}$ , and  $Y=N_B(u)=\{y_1,\ldots,y_4\}.$
By a similar argument as the proof in Claim \ref{t=14}, $X$ forms one of the Ramsey graphs $\{\Gamma_{8}, \Gamma_8', \Gamma_8''\}$, and $X\cup \{u\}$ contains two disjoint red half $K_3$'s and an independent half edge. Note that $Y$ contains either a red half  $K_3$ or two independent half  edges from Fact \ref{center-fact} (iv). If the first case is true, then we apply Lemma \ref{1-sys} with $\ell_6=3$, $\ell_7=1$ and $t=13$ to obtain $2g( {{K_t}( \omega  )} ) \le 1 - \frac{{30}}{{30t - 45\cdot 3-20}} < 7/8.$ If the second case occurs, then we apply Lemma \ref{1-sys} with $\ell_6=2$, $\ell_7=3$ and $t=13$ to obtain $2g( {{K_t}( \omega  )} ) \le 1 - \frac{{30}}{{30t - 45\cdot 2-20\cdot3}} =7/8$.
Moreover, if the equality holds, then all the other red edges in $X$ are full apart from the two independents $K_3$'s edges and these three independent half  edges mentioned above, contradicting Fact \ref{eightdagon-like-K3-K7} (i) that all red edges contained in full-neighborhood of $u$ are half.

Now, we assume that $\deg_B(u)=5$, and $\deg_R(u)=7$.  Let
$X = N_R( u ) =\{ {x_1,\ldots,x_7} \}$, and $Y= N_B( u)= \{ y_1,\ldots,y_5\}.$ From Fact \ref{center-fact} (iii), $Y$ forms a red half  $K_5$.

We first assume that $u$ has one half-neighbor in $X$, say $x_7$, and all $ux_i$ for $i\in[6]$ are full  edges. Then $X\setminus\{x_7\}$ contains a red half  $K_3$ from Fact \ref{two-ind-eg-for-K3-K7}. Then we apply Lemma \ref{1-sys} with $\ell_2=1$, $\ell_6=\ell_7=1$, and $t=13$ to obtain $2g( {{K_t}( \omega  )} ) \le 1 - \frac{{30}}{{30t - 100- 45-20}} <7/8$.

 Now we may assume that $u$ is {\bf full}-adjacent to each vertex in $X$, then from Fact \ref{two-ind-eg-for-K3-K7}, $X$ contains a {\bf red} half  $K_3$, say on $\{x_1,x_2,x_3\}$. Let $X_1=\{x_4,\ldots,x_7\}$.

 We will show $X_1$ contains at least $1$ half  edge. ~Otherwise, all edges in $X_1$ are {\bf full}. Then $X_1$ contains no red $K_3$ by Fact \ref{center-fact} (iv). Thus, from Fact \ref{T-I-E}, $X_1$ contains two independent red edges, say $x_4x_5$ and $x_6x_7$. If there exists some $i_0\in [3]$ such that $\deg_B(x_{i_0},X_1)\ge3$, then all edges contained in $N_B(x_{i_0},X_1)$ are red to avoid a blue ${K_3}$. However, then $N_B(x_{i_0},X_1)$ contains a red full  $K_3$, contradicting Fact \ref{center-fact} (iv) again. Thus for each $i\in[3]$, $\deg_B(x_i,X_1)\le2$, and so $\deg_R(x_i,X_1)\ge2$. If there is some $j_0\in [3]$ such that $\deg_R(x_{j_0},X_1)\ge3$, then the red full  edge $x_4x_5$ or $x_6x_7$ will be contained in $N_R(x_{j_0},X_1)$, which again contradicts Fact \ref{center-fact} (iv). Thus, for each $i\in[3]$, we have $\deg_R(x_i,X_1)\le2$, and so
$
\deg_R(x_i,X_1)=\deg_B(x_i,X_1)=2.
$
Moreover, from Fact \ref{center-fact} (iv), $N_R(x_i,X_1)\neq\{x_4,x_5\},\{x_6,x_7\}$  for each $i\in[3]$. By symmetry, we may assume $N_R(x_1,X_1)=\{x_4,x_7\}$, then $x_4x_7$ must be blue. Thus $N_R(x_2,X_1)\neq\{x_5,x_6\}$. Note that $N_B(x_1,X_1)=\{x_5,x_6\}$, so $x_5x_6$ must be red. This implies that $x_4x_6$ and $x_5x_7$ are blue, and so $N_R(x_2,X_1)\neq \{x_4,x_6\}, \{x_5,x_7\}$. Thus, $N_R(x_2,X_1)=\{x_4,x_7\}$. Similarly, $N_R(x_3,X_1)=\{x_4,x_7\}$. However, then $\{u,x_1,x_2,x_3,x_4\}$ would form a red $K_5$ with a full edge $ux_1$, contradicting Fact \ref{center-fact} (iii). Therefore, $X_1$ has at least $1$ half  edge as claimed.

Recall $Y$ forms a red half  $K_5$, and $\{x_1,x_2,x_3\}$ forms a red half  $K_3$, so we apply Lemma \ref{1-sys} with $\ell_2=1$, $\ell_6=\ell_7=1$, and $t=13$ to obtain $2g( {{K_t}( \omega  )} ) \le 1 - \frac{{30}}{{30t - 100- 45-20}} <7/8.$ Claim \ref{t=13} is proved.\hfill$\Box$

\smallskip
From the above claim, we have
$
 2g( {{K_t}( \omega  )} ) \le 1 - \frac{1}{t} - \frac{2}{{2t}} < 7/8
 $ by applying Lemma \ref{123} with $\lambda = 2$ and $t = 13$.
The proof for $t=13$ is complete.

\medskip\noindent
{\bf Part (III) \; $t=10,11,12$}

\medskip


We have the following claim.
\begin{claim}\label{t=10-12}
For each $u \in {K_{t}}( \omega  )$ where $t=10,11,12$, $u$ is incident to at least $1$ half  edge.
\end{claim}

In order to keep the main
line clear, we will give the proof of the above claim in the next subsection
since the proof is technical. 

For $t=12$, from Claim \ref{t=10-12} and  Lemma \ref{123}, we obtain
$
 2g( {{K_t}( \omega  )} ) \le 1 - \frac{1}{t} - \frac{2}{{2t}} = 7/8
$, where the equality holds only when  all half  edges form a perfect matching on $K_{12}(\omega)$. From Fact \ref{center-fact} (ii), if $\deg_B(x_0)=5$ for some $x_0\in V(K_{12}(\omega))$, then $N_B(x_0)$ induces a red half  $K_5$ from Fact \ref{center-fact} (iii), a contradiction. Thus, for all $x\in V(K_{12}(\omega))$,  $\deg_R(x)\ge7$, and so each vertex has at least 6 red full-neighbors, which implies that the red full-neighborhood of $x$ contains a red half  $K_3$ by Fact \ref{two-ind-eg-for-K3-K7}, a contradiction. The proof for $t=12$ is complete.
For $t=11$ and $10$, we can also apply Lemma \ref{123} and Claim \ref{t=10-12} to obtain
$ 2g( {{K_t}( \omega  )} )  < 7/8,$ respectively.

\medskip
From Parts (I)-(III), $\rho(3,7)=7/16$ as desired.
Moreover, only when $t=8$, such a {\em GBE}-graph $B_n(8)$ can achieve the Ramsey-Tur\'{a}n density $7/16$ and $\Delta(B_n(8), U(n,4,3))=o(n^2)$. \hfill$\Box$

\subsection{Proof of Claim \ref{t=10-12}}

Since the proofs of Claim \ref{t=10-12} for $t=10,11,12$ are much different, we separate Claim \ref{t=10-12} into three claims, i.e., Claim \ref{t=12}, Claim \ref{t=11}, and Claim \ref{t=10}, as follows.

\begin{claim}\label{t=12}
For each $u \in {K_{12}}( \omega  )$, $u$ is incident to at least $1$ half  edge.
\end{claim}

\noindent{\em Proof.} ~On the contrary, there exists a vertex $u$ which is {\bf full}-adjacent to all other vertices. From Fact \ref{center-fact} (ii), we have ${\deg_B}( u ) \le 5$. Suppose first ${\deg_B}( u ) \le 3$, and so ${\deg_R}( u ) \ge 8$. It follows from Fact \ref{eightdagon-like-K3-K7} (ii) that ${\deg_R}( u ) = 8$, and thus $N_R(u)$ contains a red half  $C_8$ as subgraph from Fact \ref{eightdagon-like-K3-K7} (i). Therefore, we  obtain $2g( {{K_t}( \omega  )} ) \le 1 - \frac{{30}}{{30t - 120}} =7/8$ by applying Lemma \ref{1-sys} with $\ell_1=1$ and $t=12$. Moreover, if the equality holds, then all the other red edges in $N_R{(u)}$ are full apart from this $C_8$, contradicting Fact \ref{eightdagon-like-K3-K7} (i) that all red edges contained in full-neighborhood of $u$ are half.

Next, we suppose $\deg_B(u)=5$. Let  $X=N_R(u)=\{x_1,\ldots,x_6\}$, and $Y=N_B(u)=\{y_1,\ldots,y_5\}.$
It follows from Fact \ref{two-ind-eg-for-K3-K7} that $X$ contains a red half  $K_3$ by noting $u$ is red full-adjacent to each vertex in $X$. Moreover, $Y$ forms a red half  $K_5$ from Fact \ref{center-fact} (iii). Thus, we apply Lemma \ref{1-sys} with $\ell_2=\ell_6=1$, and $t=12$ to obtain $2g( {{K_t}( \omega  )} ) \le 1 - \frac{{30}}{{30t - 100-45}} <7/8$.

It remains to check ${\deg_B}( u ) = 4$. Let $X=N_R(u)=\{x_1,\ldots,x_7\}$, and $Y=N_B(u)=\{y_1,\ldots,y_4\}.$
Note that $Y$ contains either a red half  $K_3$ or two independent half  edges from Fact \ref{center-fact} (iv), and $X$ contains a red half  $K_3$, say on $\{x_1,x_2,x_3\}$ from Fact \ref{two-ind-eg-for-K3-K7}. Let $X_1=\{x_4,\ldots,x_7\}$.

\begin{proposition}\label{t=12-1}
$X_1$ contains no red $K_3$.
\end{proposition}
\noindent{\em Proof.} ~Otherwise, $X_1$ contains a red $K_3$, which must be half from Fact \ref{center-fact} (iv) by noting $X_1\subseteq N_R(u)$ and all edges incident to $u$ are full. Recall that $Y$ contains either a red half  $K_3$ or two independent half  edges and $\{x_1,x_2,x_3\}$ is disjoint from $Y$. If $Y$ contains a red half  $K_3$, then we apply Lemma \ref{1-sys} with $\ell_6=3$, and $t=12$ to obtain $2g( {{K_t}( \omega  )} ) \le 1 - \frac{{30}}{{30t - 45\cdot 3}} < 7/8$. If $Y$ has two independent half  edges, then we apply Lemma \ref{1-sys} with $\ell_6=2$, $\ell_7=2$ and $t=12$ to obtain $2g( {{K_t}( \omega  )} ) \le 1 - \frac{{30}}{{30t - 45\cdot 2-20\cdot2}} <7/8$. \hfill$\Box$

\begin{proposition}\label{t=12-2}
$X_1$ contains $2$ independent red half  edges.
\end{proposition}
\noindent{\em Proof.} ~From Proposition \ref{t=12-1} and Fact \ref{T-I-E}, $X_1$ contains two independent {\bf red} edges, say $x_4x_5$ and $x_6x_7$. If $x_4x_5$ and  $x_6x_7$ are half  edges, we are done.

First, we suppose both $x_4x_5$ and  $x_6x_7$ are {\bf full}. Then, from Fact \ref{center-fact} (iv), $\deg_R(x_{i},X_1)\le2$ for each $i\in[3]$. If there exists some $i_0\in [3]$ with $\deg_B(x_{i_0},X_1)\ge3$, then $N_B(x_{i_0},X_1)$ forms a red $K_3$, contradicting Proposition \ref{t=12-1}. Thus, for each $i\in[3]$,
$
\deg_R(x_i,X_1)=\deg_B(x_i,X_1)=2.
$

Clearly, we have $N_R(x_i,X_1)\neq\{x_4,x_5\},\{x_6,x_7\}$  for each $i\in[3]$.  By symmetry, we may assume $N_R(x_1,X_1)=\{x_4,x_7\}$. Then $N_B(x_1,X_1)=\{x_5,x_6\}$, which implies that $x_5x_6$ is {\bf red}. Since $x_4x_5$ and $x_6x_7$ are red and $X_1$ contains no red $K_3$ from Proposition \ref{t=12-1}, we have $x_4x_6$ and $x_5x_7$ must be blue, which implies $N_R(x_2,X_1)\neq\{x_5,x_7\},\{x_4,x_6\}$ to avoid a blue $K_3$. Moreover, we have $N_R(x_2,X_1)\neq\{x_4,x_5\},\{x_6,x_7\}$ from Fact \ref{center-fact} (iv).
So
$
N_R(x_2,X_1)=\{x_4,x_7\}, \; \text{or} \;\; \{x_5,x_6\}.
$
 Similarly, $N_R(x_3,X_1)=\{x_4,x_7\}$ or $\{x_5,x_6\}$.

Recall $N_R(x_1,X_1)=\{x_4,x_7\}$. If $N_R(x_2,X_1)=N_R(x_3,X_1)=\{x_4,x_7\}$, then $\{u,x_1,x_2,x_3,x_7\}$ forms a red $K_5$ in which $ux_1$ is full, contradicting Fact \ref{center-fact} (iii). If $N_R(x_2,X_1)=N_R(x_3,X_1)=\{x_5,x_6\}$, then $\{u,x_2,x_3,x_5,x_6\}$ forms a red $K_5$ in which $ux_2$ is full, a contradiction. Thus we may assume $N_R(x_2,X_1)=\{x_5,x_6\}$, and $N_R(x_3,X_1)=\{x_4,x_7\}$. Then $N_B(x_2,X_1)=\{x_4,x_7\}$, and so $x_4x_7$ is  red to avoid a blue $K_3$. However, then $\{u,x_1,x_3,x_4,x_7\}$ forms a red $K_5$'s with a full  edge $ux_1$, again a contradiction.

It remains to check that exactly one of the {\bf red} edges $\{x_4x_5,x_6x_7\}$ is full. By symmetry, we may assume that $x_4x_5$ is {\bf full} and $x_6x_7$ is {\bf half}.
Recall $\{x_1,x_2,x_3\}$ forms a red half  $K_3$.

\begin{observation}\label{12-5}
For each $i\in[3]$, $\deg_R(x_i,X_1)=\deg_B(x_i,X_1)=2$.
\end{observation}
\noindent{\em Proof.}  To avoid a blue ${K_3}$, all edges contained in $N_B(x_{i},X_1)$ are red for $i\in[3]$, it follows from Proposition \ref{t=12-1} that for each $i\in[3]$,
\begin{align}\label{12-4}
 \deg_B(x_i,X_1)\le2,\; \text{and so}\ \deg_R(x_i,X_1)\ge2.
\end{align}

We will show the equalities in (\ref{12-4}) hold.
Otherwise, suppose $\deg_R(x_1,X_1)\ge3$ by symmetry. Then $\{x_4,x_5\}\nsubseteq N_R(x_{1},X_1)$ from Fact \ref{center-fact} (iv) since $x_4x_5$ is full. By symmetry, we may assume $\{x_4,x_6,x_7\}\subseteq N_R(x_{1},X_1)$.
Clearly, $x_{1}x_6$ and $x_{1}x_7$ are half  edges from Fact \ref{center-fact} (iv).

$x_{1}x_4$ is half; otherwise, $x_4x_6$ and $x_4x_7$ must be blue. Similarly, $x_2x_4$ and $x_3x_4$ are blue. Thus, $x_{2}x_6,x_2x_7,x_{3}x_6,x_3x_7$ are red to avoid a blue $K_3$. However, then $\{u,x_2,x_3,x_6,x_7\}$ would form a red $K_5$ with a full  edge $ux_2$, which contradicts Fact \ref{center-fact} (iii).

$x_4x_6$ is full; otherwise, $x_1x_4x_6x_7x_1$ is a half  $C_4$. Recall that $x_2x_3$ is a red half  edge, and $Y$ contains either a red half  $K_3$ or two independent red half  edges. If $Y$ contains a red half  $K_3$, then we apply Lemma \ref{1-sys} with $\ell_5=\ell_6=\ell_7=1$, and $t=12$ to obtain $2g( {{K_t}( \omega  )} ) \le 1 - \frac{{30}}{{30t - 60-45-20}} <7/8.$ If $Y$ contains two independent half  edges, then we apply Lemma \ref{1-sys} with $\ell_5=1$, $\ell_7=3$ and $t=12$ to obtain $2g( {{K_t}( \omega  )} ) \le 1 - \frac{{6}}{{6t - 12-4\cdot3}} =7/8$, and the equality does not hold since there are other half edges.

Then $x_4x_6$ must be blue; otherwise, $\{u,x_1,x_4,x_6\}$ would be a red $K_4$ in which $\{u,x_4,x_6\}$  forms a red full  $K_3$, contradicting Fact \ref{center-fact} (iv). Similarly, $x_4x_7$ is a blue full  edge.

We have $\deg_R(x_2,X_1)=\deg_R(x_3,X_1)=2$; otherwise, there would be a red $K_5$ with a full  edge, or there would be a red $K_4$ with a full $K_3$, which is impossible.

Recall that $x_4x_6$ and $x_4x_7$ are blue, so for $j\in\{2,3\}$, $N_R(x_{j},X_1)\neq\{x_5,x_7\},\{x_5,x_6\}$ to avoid a blue $K_3$. Since $\{x_6,x_7\}\subseteq N_R(x_1,X_1)$ and $x_6x_7$ is red, we have $N_R(x_{j},X_1)\neq\{x_6,x_7\}$ for $j\in\{2,3\}$ from Fact \ref{center-fact} (iii). Moreover, $N_R(x_{j},X_1)\neq\{x_4,x_5\}$ for $j\in\{2,3\}$ from Fact \ref{center-fact} (iv) as $x_4x_5$ is a red full  edge. Therefore, $N_R(x_{j},X_1)=\{x_4,x_7\}$ or $\{x_4,x_6\}$ for $j\in\{2,3\}$. However, then $\{u,x_1,x_2,x_3,x_4\}$ forms a red $K_5$ with a full  edge $ux_1$, a contradiction. Therefore, $\deg_R(x_1,X_1)\le2$. Observation \ref{12-5} is finished.\hfill$\Box$

\smallskip
Recall $\{x_1,x_2,x_3\}$ forms a red half  $K_3$, $x_4x_5$ is a red full  edge and $x_6x_7$ is a red half  edge.  Note that $\deg_R(x_i,X_1)=\deg_B(x_i,X_1)=2$ for $i\in[3]$ from Observation \ref{12-5}.
Thus for $i\in[3]$, $N_R(x_i,X_1)\neq \{x_4,x_5\}$ by Fact \ref{center-fact} (iii).

Suppose that there is a vertex in $\{x_1,x_2,x_3\}$, say $x_1$, such that $N_B(x_1,X_1)=\{x_4,x_5\}$, and so $N_R(x_1,X_1)=\{x_6,x_7\}$.
Then  for $j\in\{2,3\}$, $N_R(x_j,X_1)\neq \{x_6,x_7\}$ from Fact \ref{center-fact} (iii). Therefore, by symmetry, we may assume that $N_R(x_2,X_1)=\{x_5,x_6\}$, and so $x_4x_7$ is red to avoid a blue $K_3$. Thus $x_3x_6$ is blue; otherwise, $\{u,x_1,x_2,x_3,x_6\}$ forms a red $K_5$ with a full  edge $ux_1$, a contradiction. From Proposition \ref{t=12-1}, $x_4x_6$ and $x_5x_7$ are blue. Thus  for $j\in\{2,3\}$, $N_R(x_j,X_1)\neq \{x_5,x_7\},\{x_4,x_6\}$ to avoid a blue $K_3$.
Thus, $x_3x_4$ is red; otherwise $\{x_3,x_4,x_6\}$ would be a blue $K_3$. This implies that $x_3x_5$ must be blue; otherwise, $\{u,x_3,x_4,x_5\}$ would form a red $K_4$ with a full  $K_3$ on $\{u,x_4,x_5\}$.
 So $x_5x_6$ must be red to avoid a blue $K_3$. Moreover, $x_3x_7$ is red. We conclude that $N_R(x_2,X_1)=\{x_5,x_6\}$ and $N_R(x_3,X_1)=\{x_4,x_7\}$, and $x_5x_6$ and $x_4x_7$ are red. Thus, $x_5x_6$ and $x_4x_7$ must be half from Fact \ref{center-fact} (iv) since $u$ is red full-adjacent to $X$.

Now we assume for each $i\in[3]$, $N_B(x_i,X_1)\neq\{x_4,x_5\}$, and so $N_R(x_i,X_1)\neq\{x_6,x_7\}$. Recall $N_R(x_i,X_1)\neq\{x_4,x_5\}$. Thus we may assume $N_R(x_1,X_1)=\{x_4,x_6\}$ by symmetry, then $N_B(x_1,X_1)=\{x_5,x_7\}$, and so $x_5x_7$ is red to avoid a blue $K_3$. Since $X_1$ contains no red $K_3$ from Proposition \ref{t=12-1}, and $x_4x_5,x_6x_7$ are red, we have $x_4x_7,x_5x_6$ are blue. Thus, for $j\in[2,3]$, $N_R(x_j,X_1)\neq\{x_4,x_7\},\{x_5,x_6\}$. Therefore, for $j\in[2,3]$,
 \begin{align}\label{2,3}
 N_R(x_j,X_1)= \{x_4,x_6\}\;\text{or}\;\{x_5,x_7\}.
 \end{align}

 If one vertex in $\{x_2,x_3\}$, say $x_2$, satisfies $N_R(x_2,X_1)=\{x_5,x_7\}$, then $N_B(x_2,X_1)=\{x_4,x_6\}$, and so $x_4x_6$ must be red. Thus,  $N_R(\{x_2,x_3\},X_1)=\{x_5,x_7\}$ or $N_R(\{x_1,x_3\},X_1)=\{x_4,x_6\}$ from (\ref{2,3}). Each case would yield a red $K_5$ with a red full edge, a contradiction. Therefore, for $j\in[2,3]$, there must be $N_R(x_j,X_1)=\{x_4,x_6\}$. However, then $\{u,x_1,x_2,x_3,x_4\}$ would be a red with a red full edge $ux_1$,  again a contradiction.
Proposition \ref{t=12-2} is complete.\hfill$\Box$

\medskip\noindent
{\bf Continue the proof of Claim \ref{t=12}.} Recall $Y$ contains either a red half  $K_3$ or two independent half  edges, and $\{x_1,x_2,x_3\}$ forms a red half  $K_3$, and $X_1$ contains two independent red half  edges. If the former is true, then we apply Lemma \ref{1-sys} with $\ell_6=2$, $\ell_7=2$ and $t=12$ to obtain
$2g( {{K_t}( \omega  )} ) \le 1 - \frac{{30}}{{30t - 45\cdot2-20\cdot2}} <7/8.$ If the later occur, then we apply with Lemma \ref{1-sys} with $\ell_6=1$, $\ell_7=4$,  and $t=12$ to obtain $2g( {{K_t}( \omega  )} ) \le 1 - \frac{{30}}{{30t - 45-20\cdot4}} <7/8$.
Claim \ref{t=12} is proved.\hfill$\Box$

\begin{claim}\label{t=11}
For each $u \in {K_{11}}( \omega  )$, $u$ is incident to at least $1$ half  edge.
\end{claim}
\noindent{\em Proof.} ~On the contrary, there is a vertex $u$ which is {\bf full}-adjacent to all other vertices. From Fact \ref{center-fact} (ii), we have ${\deg_B}( u ) \le 5$.
If ${\deg_B}( u ) \le 2$, then by a similar argument as that of Claim \ref{t=12}, we obtain $2g( {{K_t}( \omega  )} )  <7/8$.
   If ${\deg_R}( u ) ={\deg_B}( u ) =5$, then from Fact \ref{center-fact} (iii), $Y$ forms a red half  $K_5$.  Thus, we apply Lemma \ref{1-sys} with $\ell_2=1$, and $t=11$ to get   $2g( {{K_t}( \omega  )} ) \le 1 - \frac{{30}}{{30t - 100}} <7/8$.

Now suppose ${\deg_R}( u ) = 6$. Let $X=N_R(u)=\{x_1,\ldots,x_6\}$, and $Y=N_B(u)=\{y_1,\ldots,y_4\}.$
Note that $Y$ contains either a red half  $K_3$ or two independent half  edges from Fact \ref{center-fact} (iv) and $X$ contains a red half  $K_3$, say on $\{x_1,x_2,x_3\}$, from Fact \ref{two-ind-eg-for-K3-K7}. Let $X_1=\{x_4,x_5,x_6\}$.

We claim that $X$ contains either a red half  $C_4$ or two disjoint red half  $K_3$'s.
 Suppose that $X$ contains no red half  $C_4$, then $\deg_R(x,X\setminus{X_1})\le 1$ for each $x\in X_1$; otherwise, there would exist a red $K_4$ with a full  $K_3$, a contradiction. Thus for each $x\in X_1$, $\deg_B(x,X\setminus{X_1})\ge 2$, and so $N_B(\{x',x''\},X\setminus{X_1})\neq\emptyset$ for $x',x''\in X_1$, which implies that $X_1$ induces a red $K_3$ to avoid a blue $K_3$. Since $u$ is red full-adjacent to each vertex in $X$, by Fact \ref{center-fact} (iv), $X_1$ forms a red half  $K_3$. This together with $X\setminus{X_1}$ yield two disjoint red half  $K_3$'s.

From the above, if $Y$ contains a red half  $K_3$, then we apply Lemma \ref{1-sys} with $\ell_5=\ell_6=1$ and $t=11$ to get $2g( {{K_t}( \omega  )} ) \le 1 - \frac{{6}}{{6t - 12-9}} <7/8$, or with $\ell_6=3$ and $t=11$ to obtain $2g( {{K_t}( \omega  )} ) \le 1 - \frac{{6}}{{6t - 9\cdot3}} <7/8$; if $Y$ contains two independent half  edges, we apply Lemma \ref{1-sys} with $\ell_5=1$, $\ell_7=2$ and $t=11$ to obtain $2g( {{K_t}( \omega  )} ) \le 1 - \frac{{6}}{{6t - 12-4\cdot2}} <7/8$, or with $\ell_6=\ell_7=2$ and $t=11$ to obtain $2g( {{K_t}( \omega  )} ) \le 1 - \frac{{6}}{{6t - 9\cdot2-4\cdot2}} <7/8$.

Finally, we assume ${\deg_R}( u ) =7$. Let $X=N_R(u)=\{x_1,\ldots,x_7\}.$
From Fact \ref{two-ind-eg-for-K3-K7}, $X$ contains a \textbf{red half}-$K_3$, say on $\{x_1,x_2,x_3\}$. 
Let $X_1=\{x_4,\ldots,x_7\}$.

We claim that $X_1$ contains no red $K_3$. Otherwise, there is a red $K_3$ in $X_1$, say on $T=\{x_4,x_5,x_6\}$, which must be half from Fact \ref{center-fact} (iv) by noting $X_1\subseteq N_R(u)$. Thus, we apply Lemma \ref{1-sys} with $\ell_6=2$ and $t=11$ to obtain $2g( {{K_t}( \omega  )} ) \le 1 - \frac{{30}}{{30t - 45\cdot 2}} =7/8$, where the equality holds only when all edges but that on $\{x_1,x_2,x_3\}$ and $T$ are full. Together with Fact \ref{center-fact} (iv), we have for $i\in[3]\cup \{7\}$, $\deg_R(x_i,T)\le 1$, and so $\deg_B(x_i,T)\ge 2$. Thus, $N_B(\{x_j,x_7\},T)\neq \emptyset$ for each $j\in [3]$, which implies that $x_jx_7$ is red to avoid a blue $K_3$. Thus $\{u,x_1,x_2,x_3,x_7\}$ would be a red $K_5$ with a full  edge $ux_1$, contradicting Fact \ref{center-fact} (iii).

Therefore, from Fact \ref{T-I-E}, $X_1$ contains two independent \textbf{red} edges, say $x_4x_5$ and $x_6x_7$. Since $X_1$ contains no red $K_3$, there are at most two red edges in $X_1$ beside $x_4x_5$ and $x_6x_7$. We have three cases.

\medskip

\textbf{Case 1:} Beside $x_4x_5$ and $x_6x_7$, all other edges in $X_1$ are blue.

\medskip
For each $i\in[3]$, $\deg_B(x_i,X_1)\le 2$ to avoid a blue $K_3$, and so $\deg_R(x_i,X_1)\ge 2$. We claim that $\{x_4,x_5\}$ or $\{x_6,x_7\}$ is contained in $N_R(x_i)$  for each $i\in[3]$. It is clear if $\deg_R(x_i,X_1)\ge 3$, so we assume $\deg_R(x_i,X_1)= 2$.
Note that all edges but $x_4x_5$ and $x_6x_7$ in $X_1$ are blue, so $N_R(x_i,X_1)=\{x_4,x_5\}$ or $\{x_6,x_7\}$ to avoid a blue $K_3$. Thus, there must exist at least $2$ vertices in $\{x_1,x_2,x_3\}$, say $\{x_1,x_2\}$, such that $\{x_4,x_5\}\subseteq N_R(\{x_1,x_2\},X_1)$ or $\{x_6,x_7\}\subseteq N_R(\{x_1,x_2\},X_1)$, contradicting Fact \ref{center-fact} (iii) by noting $u$ is red full-adjacent to $X$.

\medskip

\textbf{Case 2:} Beside $x_4x_5$ and $x_6x_7$, there exists exactly one more red edge in $X_1$, say $x_4x_7$.

\medskip
Then $x_4x_6, x_5x_6, x_5x_7$ are \textbf{blue}. Since $K_t(w)$ has no blue $K_3$, we have for $i\in[3]$,
\begin{align}\label{11-deg}
 \deg_B(x_i,X_1)\le 2, \;\;\text{and} \;\; \deg_R(x_i,X_1)\ge 2.
\end{align}

Clearly, $\deg_R(x_i,X_1)\le3$ for each $i\in[3]$; otherwise, assume $\deg_R(x_1,X_1)=4$ by symmetry. Since $x_4x_5,x_6x_7,x_4x_7,x_1x_2$ are red, we have $\deg_R(x_2,X_1)=2$ from Fact \ref{center-fact} (iii). Moreover, $N_R(x_2,X_1)\neq\{x_4,x_5\},\{x_6,x_7\},\{x_4,x_7\}$. Since $x_4x_6, x_5x_6, x_5x_7$ are blue, we have $N_R(x_2,X_1)\neq\{x_5,x_7\},\{x_4,x_7\},\{x_4,x_6\}$ to avoid a blue $K_3$.
 Thus, $N_R(x_2,X_1)=\{x_5,x_6\}$. Similarly, $N_R(x_3,X_1)=\{x_5,x_6\}$. However, then $\{u,x_1,x_2,x_3,x_5\}$ would form a red $K_5$ with a red full  edge $ux_1$, a contradiction.

\begin{observation}\label{t=11-3}
At most one vertex in $\{x_1,x_2,x_3\}$ has three red neighbors in $X_1$.
\end{observation}
\noindent{\em Proof.} On the contrary, we may assume $\deg_R(x_1,X_1)=\deg_R(x_2,X_1)=3$ by symmetry.

If $N_R(x_1,X_1)=\{x_4,x_5,x_6\}$, then $N_R(x_2,X_1)=\{x_4,x_6,x_7\}$ or $N_R(x_2,X_1)=\{x_5,x_6,x_7\}$; otherwise, $\{x_4,x_5\}\subseteq N_R(\{x_1,x_2\},X_1)$, then there would be a red $K_5$ with a full edge by noting $x_4x_5,x_1x_2$ are red. Thus $\deg_R(x_3,X_1)=2$ from (\ref{11-deg}); otherwise, $\{x_4,x_5\}\subseteq N_R(\{x_1,x_3\},X_1)$ or $\{x_6,x_7\}\subseteq N_R(\{x_2,x_3\},X_1)$, contradicting Fact \ref{center-fact} (iii) again. Thus, $N_R(x_3,X_1)\neq\{x_4,x_5\}$, $\{x_6,x_7\}.$
Moreover, since $x_4x_6,x_5x_6,x_5x_7$ are blue, we have $N_R(x_3,X_1)\neq\{x_5,x_7\}$, $\{x_4,x_7\}$, $\{x_4,x_6\}$ to avoid a blue $K_3$. Therefore, $N_R(x_3,X_1)=\{x_5,x_6\}$. However, then $\{u,x_1,x_2,x_3,x_6\}$ would form a red $K_5$ with a full  edge $ux_1$, a contradiction.

If $N_R(x_1,X_1)=\{x_4,x_5,x_7\}$, then $N_R(x_2,X_1)=\{x_5,x_6,x_7\}$; otherwise, we have that $\{x_4,x_5\}\subseteq N_R(\{x_1,x_2\},X_1)$ or $\{x_4,x_7\}\subseteq N_R(\{x_1,x_2\},X_1)$, contradicting Fact \ref{center-fact} (iii) since $x_1x_2$, $x_4x_5$, $x_4x_7$ are red. A similar argument as above yields $N_R(x_3,X_1)=\{x_5,x_6\}$. However, then $\{u,x_1,x_2,x_3,x_5\}$ would form a red $K_5$ with a full  edge $ux_1$, a contradiction.

All other cases, i.e., $N_R(x_1,X_1)=\{x_5,x_6,x_7\}$ and $N_R(x_1,X_1)=\{x_4,x_6,x_7\}$,  are also impossible, and the proofs are similar. \hfill$\Box$

\smallskip

From Observation \ref{t=11-3}, we may first assume $\deg_R(x_1,X_1)=3$, and $\deg_R(x_j,X_1)=2$ for $j\in\{2,3\}$ by symmetry. Since $x_4x_6,x_5x_6,x_5x_7$ are blue,  $N_R(x_j,X_1)\neq\{x_5,x_7\},\{x_4,x_7\},\{x_4,x_6\}$ for $j\in\{2,3\}$.
If $N_R(x_1,X_1)=\{x_4,x_5,x_6\}$,  then $N_R(x_j,X_1)\neq\{x_4,x_5\}$ for $j\in\{2,3\}$ from Fact \ref{center-fact} (iii), implying $N_R(x_j,X_1)=\{x_5,x_6\}$ or $\{x_6,x_7\}$. However, then $\{u,x_1,x_2,x_3,x_6\}$ would form a red $K_5$ with full  edge $ux_1$, a contradiction. By symmetry, it remains to check that $N_R(x_1,X_1)=\{x_4,x_5,x_7\}$. Then $\{x_1,x_4,x_5\}$ and $\{x_1,x_4,x_7\}$ would be two red half  $K_3$'s from Fact \ref{center-fact} (iv) since $u$ is red full-adjacent to each vertex in $X$, which implies that $x_1x_5x_4x_7x_1$ forms a red half  $C_4$. Similarly, we have $N_R(x_j,X_1)=\{x_5,x_6\}$ or $\{x_6,x_7\}$ for $j\in\{2,3\}$, which implies that $\{x_2,x_3,x_6\}$ forms a red half  $K_3$ from Fact \ref{center-fact} (iv). Thus, we apply Lemma \ref{1-sys} with $\ell_5=\ell_6=1$, and $t=11$ to obtain
$2g( {{K_t}( \omega  )} ) \le 1 - \frac{{6}}{{6t - 12-9}} <7/8.$

Now suppose $\deg_R(x_i,X_1)=2$ for each $i\in[3]$. Since $x_4x_6,x_5x_6,x_5x_7$ are blue, we have for $i\in[3]$, $N_R(x_i,X_1)\neq\{x_5,x_7\},\{x_4,x_7\},\{x_4,x_6\}$ to avoid a blue $K_3$. Thus,
$N_R(x_i,X_1)=\{x_4,x_5\}$, $\{x_5,x_6\}$, or $\{x_6,x_7\}.$
If there are distinct $i,j\in[3]$ such that $N_R(x_i,X_1)=\{x_4,x_5\}$ and $N_R(x_j,X_1)=\{x_6,x_7\}$, then $\{x_i,x_4,x_5\}$ and $\{x_j,x_6,x_7\}$ form two disjoint red half  $K_3$'s from Fact \ref{center-fact} (iv) by noting $u$ is red full-adjacent to each vertex in $X$. Thus, we apply Lemma \ref{1-sys} with $\ell_6=2$, and $t=11$ to obtain
$2g( {{K_t}( \omega  )} ) \le 1 - \frac{{30}}{{30t - 45\cdot 2}} =7/8$, where the equality does not hold
since $x_1x_2$ and $x_1x_3$ are half.
Therefore, for all $i\in[3]$, we have that $N_R(x_i,X_1)\in\{\{x_4,x_5\},\{x_5,x_6\}\}$ or $N_R(x_i,X_1)\in\{\{x_5,x_6\},\{x_6,x_7\}\}$. Thus, $x_5\in N_R(\{x_1,x_2,x_3\},X_1)$ or $x_6\in N_R(\{x_1,x_2,x_3\},X_1)$ implying that either $\{u,x_1,x_2,x_3,x_5\}$ or $\{u,x_1,x_2,x_3,x_6\}$ forms a red $K_5$'s with a full  edge $ux_1$, a contradiction.

\medskip

\textbf{Case 3:} Beside $x_4x_5$ and $x_6x_7$, there are another two red edges in $X_1$.

\medskip

Recall $X_1$ contains no red $K_3$. Then $X_1$ contains a red $C_4$, we may assume that $x_4x_5x_6x_7x_4$ is the {\bf red} $C_4$ by symmetry, and so $x_4x_6$ and $x_5x_7$ are {\bf blue}. Since there is no blue $K_3$, we have for $i\in[3]$, $\deg_B(x_i,X_1)\le 2$, and so
$
\deg_R(x_i,X_1)\ge 2.
$

For each $i\in[3]$, $\deg_R(x_i,X_1)\le 3$. Otherwise, suppose $\deg_R(x_1,X_1)=4$ by symmetry. Since $x_4x_5x_6x_7x_4$ is a red $C_4$ and $x_1x_2$ is red, from Fact \ref{center-fact} (iii), $\deg_R(x_2,X_1)=2$, and $N_R(x_2,X_1)\neq\{x_4,x_5\},\{x_5,x_6\},\{x_6,x_7\},\{x_4,x_7\}$. Moreover,  since $x_4x_6,x_5x_7$ are blue, $N_R(x_2,X_1)\neq\{x_5,x_7\},\{x_4,x_6\}$. A contradiction.

\begin{observation}\label{t=11-4}
At most one vertex in $\{x_1,x_2,x_3\}$ has three red neighbors in $X_1$.
\end{observation}
\noindent{\em Proof.} ~Otherwise, suppose $\deg_R(x_1,X_1)=\deg_R(x_2,X_1)=3$ by symmetry. Since $x_1x_2$ is red, and $x_4x_5x_6x_7x_4$ is a red $C_4$, we have that $N_R(\{x_1,x_2\},X_1)$ contains no red edge by Fact \ref{center-fact} (iii). Thus we may assume  $N_R(x_1,X_1)=\{x_4,x_5,x_6\}$ and $N_R(x_2,X_1)=\{x_4,x_6,x_7\}$ by symmetry. Again by Fact \ref{center-fact} (iii), $\deg_R(x_3,X_1)=2$ and so $N_R(x_3,X_1)\in \{\{x_4,x_6\},\{x_5,x_7\}\}$ since $\{x_1x_3,x_2x_3\}$ are red, which is impossible since $x_4x_6$ and $x_5x_7$ are blue.\hfill$\Box$

\medskip

From Observation \ref{t=11-4}, we first assume that one vertex in $\{x_1,x_2,x_3\}$ has three red neighbors in $X_1$. By symmetry, we may assume $N_R(x_1,X_1)=\{x_4,x_5,x_7\}$, then $\{x_1,x_4,x_5\}$ and $\{x_1,x_4,x_7\}$ form two red half  $K_3$'s by Fact \ref{center-fact} (iv) since  $u$ is red full-adjacent to each vertex in $X$. Thus $x_1x_5x_4x_7x_1$ is a red half  $C_4$. For $j\in \{2,3\}$, since $\deg_R(x_j,X_1)=2$ and $\{x_4x_6,x_5x_7\}$ are blue, $N_R(x_j,X_1)\neq\{x_5,x_7\},\{x_4,x_6\}$. Moreover, $N_R(x_j,X_1)\neq\{x_4,x_5\},\{x_4,x_7\}$ for $j\in\{2,3\}$ by noting $x_1x_j$ are red and Fact \ref{center-fact} (iii). Thus, $N_R(x_j,X_1)=\{x_5,x_6\}$ or $\{x_6,x_7\}$ for $j\in \{2,3\}$, and so $\{x_2,x_3,x_6\}$ would be a red half  $K_3$. Therefore, we apply Lemma \ref{1-sys} with $\ell_5=\ell_6=1$,  and $t=11$ to obtain $2g( {{K_t}( \omega  )} ) \le 1 - \frac{{6}}{{6t - 12-9}} <7/8.$

We now assume $\deg_R(x_i,X_1)=2$ for $i\in[3]$. Since $x_4x_6$ and $x_5x_7$ are blue, we obtain that  for $i\in[3]$, $N_R(x_i,X_1)\neq\{x_5,x_7\},\{x_4,x_6\}$ to avoid a blue $K_3$, and so
$N_R(x_i,X_1)=\{x_4,x_5\}$, $\{x_4,x_7\}$, $\{x_5,x_6\}$, or $\{x_6,x_7\}.$
If there are $x',x''\in\{x_1,x_2,x_3\}$ with $N_R(\{x',x''\},X_1)=\emptyset$, then $\{x'\}\cup N_R(x',X_1)$ and $\{x''\}\cup N_R(x'',X_1)$ would form two disjoint red half  $K_3$'s from Fact \ref{center-fact} (iv) since $u$ is red full-adjacent to each vertex in $X$. Thus, we apply Lemma \ref{1-sys} with $\ell_6=2$, and $t=11$ to obtain $2g( {{K_t}( \omega  )} ) \le 1 - \frac{{30}}{{30t - 45\cdot 2}} =7/8,$ where the equality does not hold since $x_1x_2$ and $x_1x_3$ are half. Therefore, for distinct $x',x''\in\{x_1,x_2,x_3\}$, $N_R(\{x',x''\},X_1)\neq\emptyset$. By symmetry, we may assume that $N_R(x_1,X_1)=\{x_4,x_5\}$, $N_R(x_2,X_1)=\{x_5,x_6\}$. Then we must have $N_R(x_3,X_1)=\{x_4,x_5\}$ or $\{x_5,x_6\}$, yielding a red $K_5$ with a full edge, a contradiction.

The final case when $\deg_R(u)=7$ is finished, and so Claim \ref{t=11} is complete.\hfill$\Box$

\begin{claim}\label{t=10}
For each $u \in {K_{10}}( \omega  )$, $u$ is incident to at least $1$ half  edge.
\end{claim}

Let us begin with the following fact.

\begin{fact}\label{t=10-fact}
(i) $K_{10}(w)$ contains no half  $(K_3\sqcup K_2)$.

(ii) If $K_{10}(w)$ contains three independent half  edges, then all the other edges are full.

(iii) If $K_{10}(w)$ contains a half  $C_4$, then all the other edges are full.

\end{fact}
\noindent{\em Proof.} ~(i) If $K_{10}(w)$ contains a half  $(K_3\sqcup K_2)$, then we have $2g( {{K_t}( \omega  )} ) \le 1 - \frac{{6}}{{6t - 9-4}} <7/8$ by applying Lemma \ref{1-sys} with $\ell_6=\ell_7=1$ and $t=10$, a contradiction.

(ii) Let $x_1x_2,x_3x_4,x_5x_6$ be the $3$ independent half edges and $x_ix_j$ be another half edge in $K_{10}(\omega)$, for $i,j\in[10]$.
If $|\{1,2,3,4,5,6\}\cap\{i,j\}|=0$, then we apply Lemma \ref{1-sys} with $\ell_7=4$ and $t=10$ to obtain $2g( {{K_t}( \omega  )} ) \le 1 - \frac{{6}}{{6t - 4\cdot4}} <7/8$. If $|\{1,2,3,4,5,6\}\cap\{i,j\}|=1$, by symmetry, we assume $i=1,j=7$. Moreover, we may assume that all the other edges are full.
Thus, the corresponding matrix of $K_{10}(w)$ is

\begin{center}
$A=\begin{pmatrix}
0 & 1/2 & 1 & 1 & 1 & 1 & 1/2 & 1 & 1 & 1\\
1/2 & 0 & 1 & 1 & 1 & 1 & 1 & 1 & 1 & 1\\
1 & 1 & 0 & 1/2 & 1 & 1 & 1 & 1 & 1 & 1\\
1 & 1 & 1/2 & 0 & 1 & 1 & 1 & 1 & 1 & 1\\
1 & 1 & 1 & 1 & 0 & 1/2 & 1 & 1 & 1 & 1\\
1 & 1 & 1 & 1 & 1/2 & 0 & 1 & 1 & 1 & 1\\
1/2 & 1 & 1 & 1 & 1 & 1 & 0 & 1 & 1 & 1\\
1 & 1 & 1 & 1 & 1 & 1 & 1 & 0 & 1 & 1\\
1 & 1 & 1 & 1 & 1 & 1 & 1 & 1 & 0 & 1\\
1 & 1 & 1 & 1 & 1 & 1 & 1 & 1 & 1 & 0
\end{pmatrix}$.
\end{center}
Let  $\textbf{u}=(u_1,\ldots,u_{10})$ ($u_i\ge0$ and $\sum_{i=1}^{10}u_i=1)$. Then a standard computation yields that $2g(K_{t}(\omega))=\max \{\textbf{u}A\textbf{u}^{T}\}\approx 0.8696<7/8$.
Similarly, if $|\{i,j\}\cap\{1,2,3,4\}|=2$, then we have $2g(K_{t}(\omega))=\max \{\textbf{u}A\textbf{u}^{T}\}\approx0.8707<7/8$.

(iii) Suppose that there exist a half  $C_4$ and some other half edges. Let $x_1x_2x_3x_4x_1$ be the half  $C_4$ and $x_ix_j$ be another half  edge in $K_{10}(\omega)$, for some $i,j\in[10]$. If $|\{i,j\}\cap\{1,2,3,4\}|=0$, i.e., $ x_ix_j$ is disjoint from $C_4$, then we apply Lemma \ref{1-sys} with $\ell_5=\ell_7=1$, and $t=10$ to obtain  $2g( {{K_t}( \omega  )} ) \le 1 - \frac{{6}}{{6t - 12-4}} <7/8.$
If $|\{i,j\}\cap\{1,2,3,4\}|=1$, then  we may assume $i=1$ and $j=5$ by symmetry. Note that the density will not decrease if we replace a half edge by a full edge, so we may assume that all the other edges are full. Let $A$ be the corresponding weighted adjacent matrix of $K_{10}(\omega)$, and
let $\textbf{u}=(u_1,\ldots,u_{10})$ ($u_i\ge0$ and $\sum_{i=1}^{10}u_i=1)$. Then $\textbf{u}A\textbf{u}^{T}$ takes the maximum $7/8$ only when $\textbf{u}=(0,1/8,0,1/8,1/8,1/8,1/8,1/8,1/8,1/8)$. Note that the distribution of $K_{t}(\omega)$ satisfies ${u_i} > 0$ for $i\in[t]$ since $B_n(t)$ is dense. Thus, $2g(K_{t}(\omega))=\max \{\textbf{u}A\textbf{u}^{T}\}<7/8$ as desired. If $|\{i,j\}\cap\{1,2,3,4\}|=2$, then  we can similarly obtain $2g(K_{t}(\omega))=\max \{\textbf{u}A\textbf{u}^{T}\}<7/8$.\hfill$\Box$

\medskip
Now we give the proof for Claim \ref{t=10}.

\medskip
\noindent{\em Proof of Claim \ref{t=10}.} ~On the contrary, there exists a vertex $u$ which is {\bf full}-adjacent to all other vertices. From Fact \ref{center-fact} (ii), we have ${\deg_B}( u ) \le 5$.
If ${\deg_B}( u ) \le 1$, then by a similar argument as that of Claim \ref{t=12}, we obtain $2g( {{K_t}( \omega  )} )  <7/8$.
If $\deg_B(u)=2$, then $\deg_R(u)=7$. For this case, we consider the red-neighborhood of $u$. By a similar argument as that of Claim \ref{t=11}, we obtain $2g( {{K_t}( \omega  )} ) <7/8$.
For the case when ${\deg_B}( u ) = 5$, we let $Y=N_B(u)$. By Fact \ref{center-fact} (iii), $Y$ induces a red half  $K_5$.  Thus, we again apply Lemma \ref{1-sys} to obtain $2g( {{K_t}( \omega  )} ) \le 1 - \frac{{30}}{{30t - 100}} <7/8$.

Now suppose $\deg_B(u)=3$, and so $\deg_R(u)=6$. Let $X=N_R(u)=\{x_1,\ldots,x_6\}$. A similar argument as in Claim \ref{t=11} yields that $X$ contains either a red half  $C_4$ or two disjoint red half  $K_3$'s.
If $X$ contains two disjoint red half  $K_3$, then we apply Lemma \ref{1-sys} with $\ell_6=2$ and $t=10$ to obtain $2g( {{K_t}( \omega  )} ) \le 1 - \frac{{6}}{{6t - 9\cdot2}} <7/8$. If $X$ contains a red half  $C_4$, say $x_1x_2x_3x_4x_1$, then we apply Lemma \ref{1-sys} with $\ell_5=1$, and $t=10$ to obtain $2g( {{K_t}( \omega  )} ) \le 1 - \frac{{6}}{{6t - 12}} =7/8.$ Suppose that the equality holds. Then all edges but that of $C_4$ are full, and so $x_1x_3,x_2x_4$ must be blue; otherwise, there would be a red $K_4$ with a full  $K_3$, a contradiction. Since there is no blue $K_3$, we have $\deg_R(x_5,\{x_1,\dots,x_4\})\ge 2$ and $x_5$ is red-adjacent to two consecutive vertices in $C_4$. However, then there would be a red $K_4$ with a full  $K_3$, again a contradiction.

It remains to check the case when ${\deg_B}( u ) = 4$. Let $X=N_R(u)=\{x_1,\ldots,x_5\}$, and $ Y=N_B(u)=\{y_1,\ldots,y_4\}.$
Then $Y$ induces a red $K_4$. From Fact \ref{center-fact} (iv), $Y$ contains either a red half  $K_3$ or two independent half  edges.

\begin{fact}\label{t=10-inside fact}
Suppose that $X$ induces a pentagonlike coloring with the red cycle $C_5$. Then

(i) for each $y\in Y$, $\deg_B(y,X)\le 2$.

(ii) $N_R(y,X)$ are consecutive vertices in the $C_5$.
\end{fact}

In the following, all the summations of the subscripts are taken {\bf modular} ${\bf 5}$.
\begin{proposition}\label{t=10-1}
$X$ induces a pentagonlike coloring, and all red edges in $X$ are full.
\end{proposition}
\noindent{\em Proof.} ~We first show $X$ induces a pentagonlike coloring.  Since there is no blue $K_3$, from Lemma \ref{jian}, it suffices to show that $X$ contains no red $K_3$. Otherwise,  all edges of this $K_3$ must be half from Fact \ref{center-fact} (iv). Moreover, there exists a half edge in $Y$. This contradicts Fact \ref{t=10-fact} (i). Thus $X$ induces a pentagonlike coloring. Let $x_1x_2\cdots x_5x_1$ be the red $C_5$.

Suppose that $X$ contains a red half  edge, say $x_1x_2$. Then, from Fact \ref{t=10-fact} (i),  $Y$ contains two independent red half edges, say $y_1y_2,y_3y_4$.
Thus all the other edges are full from Fact \ref{t=10-fact} (ii).

From Fact \ref{t=10-inside fact}, there exist $i_1,i_2\in [4]$ such that $N_R(\{y_{i_1},y_{i_2}\},X)$ contains a red edge $x_{j_1}x_{j_1+1}$, where $j_1\in[5]$. Then $\{x_{j_1},x_{j_1+1},y_{i_1},y_{i_2}\}$ forms a red $K_4$, and $x_{j_1}x_{j_1+1}$ and $y_{i_1}y_{i_2}$ must be half to avoid a red $K_4$ with a full $K_3$. Thus, $x_{j_1}x_{j_1+1}=x_1x_2$ and $y_{i_1}y_{i_2}=y_1y_2$ or $y_3y_4$. Assume $y_{i_1}y_{i_2}=y_1y_2$ by symmetry. Since all edges but $x_1x_2,y_1y_2,y_3y_4$ are full, $N_R(\{x_1,x_2\},Y)=\{y_1,y_2\}$ from Fact \ref{center-fact} (iii). Moveover, $N_R(\{x_i,x_{i+1}\},Y)\le1$ for each $i\in[2,5]$ to avoid a red $K_4$ with a full  $K_3$. Note that for $i\in[4]$, $\deg_R(y_i,X)\ge 3$ and $N_R(y_i,X)$ consists of consecutive vertices in $C_5$ from Fact \ref{t=10-inside fact}. Thus, there exists $y_{j_1}y_{j_2}$ with $\{j_1,j_2\}\neq\{1,2\}$ such that $N_R(\{y_{j_1},y_{j_2}\},X)$ contains a red full edge, yielding a red $K_4$ with a full  $K_3$, a contradiction again. \hfill$\Box$

\smallskip

From Proposition \ref{t=10-1}, let $x_1x_2\cdots x_5x_1$ be the \textbf{red full  $C_5$}. 
From Fact \ref{t=10-inside fact}, there exist $i_1,i_2\in [4]$ such that $N_R(\{y_{i_1},y_{i_2}\},X)$ contains a red edge. We may assume that $\{x_1,x_2\}\subseteq N_R(\{y_1,y_2\},X)$ by symmetry. Then
  \begin{align}\label{t=10,deg-x_1x_2}
N_R(\{x_1,x_2\},Y)=\{y_1,y_2\},
\end{align}
since otherwise there would be a red $K_5$ with a full edge $x_1x_2$, a contradiction.

\medskip

\textbf{Case 1:} $y_1y_2$ is full.

\medskip

Since $\{x_1,x_2,y_1,y_2\}$ forms a red $K_4$, and $x_1x_2,y_1y_2$ are red full  edges, then there also exist two independent half  edges from Fact \ref{center-fact} (iv), say $x_1y_1$ and $x_2y_2$ by symmetry.

\begin{proposition}\label{t=10-2}
$y_3y_4$ is full.
\end{proposition}
\noindent{\em Proof.} Otherwise, $x_1y_1,x_2y_2,y_3y_4$ are independent half  edges, then all the other edges are full from Fact \ref{t=10-fact} (ii).
Thus, $N_R(\{x_i,x_{i+1}\},Y)\le1$ for $i\in[2,5]$ from Fact \ref{center-fact} (iv) by noting $x_1x_2\cdots x_5x_1$ is a red full  $C_5$. Note that for $i\in[4]$, $\deg_R(y_i,X)\ge 3$ and $N_R(y_i,X)$ consists of consecutive vertices in $C_5$ from Fact \ref{t=10-inside fact}, and $N_R(\{x_1,x_2\},Y)=\{y_1,y_2\}$ from (\ref{t=10,deg-x_1x_2}), so there must exist $y_{j_1},y_{j_2}$ with $\{j_1,j_2\}\neq\{1,2\}$ such that $N_R(\{y_{j_1},y_{j_2}\},X)$ contains a red edge, yielding a red $K_4$ with a full  $K_3$, a contradiction.\hfill$\Box$

\begin{proposition}\label{t=10-3} Relabel $y_3$ and $y_4$ if necessary,  we have that
$N_B(y_3,X)=\{x_2,x_3\}$, and $N_B(y_4,X)=\{x_1,x_5\}$.
\end{proposition}
\noindent{\em Proof.} $\{x_3,x_4\}\nsubseteq N_R(\{y_3,y_4\},X)$; otherwise, $\{x_3,x_4,y_3,y_4\}$ must contain two independent red half  edges from Fact \ref{center-fact} (iv) since $x_3x_4$ and $y_3y_4$ are red full edges. This contradicts Fact \ref{t=10-fact} (ii) since $x_1y_1$ and $x_2y_2$ are red half edges. Similarly, $\{x_4,x_5\}\nsubseteq N_R(\{y_3,y_4\},X)$.

$\{x_1,x_2\}\nsubseteq N_R(y_3,X)$; otherwise $\{x_1,x_2,y_1,y_2,y_3\}$ would be a red $K_5$ with a full edge.
Similarly, $ \{x_1,x_2\}\nsubseteq N_R(y_4,X).$
If $\deg_R(y_3,X)\ge 4$, then $\{x_3,x_4\}\subseteq N_R(\{y_3,y_4\},X)$ or $\{x_4,x_5\}\subseteq N_R(\{y_3,y_4\},X)$ from Fact \ref{t=10-inside fact} that $\deg_R(y_j,X)\ge3$ for $j\in[4]$, a contradiction. Thus $\deg_R(y_3,X)=3$, and similarly $\deg_R(y_4,X)=3$. By symmetry, we may assume that $N_R(y_3,X)=\{x_1,x_5,x_4\}$ and $N_R(y_4,X)=\{x_2,x_3,x_4\}$ as desired. \hfill$\Box$

\smallskip

Since $Y$ forms a red $K_4$, and $y_1y_2,y_3y_4$ are red full  edges, there also exist two independent half  edges in $Y$ from Fact \ref{center-fact} (iv).

Suppose first $y_1y_4$ and $y_2y_3$ are half. It follows from Fact \ref{t=10-inside fact} and (\ref{t=10,deg-x_1x_2}) that $x_3\in N_R(y_1,X)$ or $x_5\in N_R(y_1,X)$. If $x_3\in N_R(y_1,X)$, then $x_3y_4$ or $x_2y_4$ is half; otherwise,  from Proposition \ref{t=10-3}, $\{x_2,x_3,y_1,y_4\}$ would be a red $K_4$ with a red full  $K_3$, a contradiction. Note that $x_1y_1,x_2y_2,y_2y_3$ are half edges, so either $x_1y_1,y_2y_3,x_3y_4$ form $3$ independent half edges  or $x_1y_1,y_2y_3,x_2y_4$ form $3$ independent half edges, contradicting Fact \ref{t=10-fact} (ii) by noting $y_1y_4$ is also half. Similarly, $x_5\in N_R(y_1,X)$ is impossible.

Now suppose $y_1y_3$ and $y_2y_4$ are half. This will lead to a contradiction similar as above by noting $x_2y_4$ is a red full edge from Proposition \ref{t=10-3}.

\medskip

\textbf{Case 2:} $y_1y_2$ is half.

\begin{proposition}\label{t=10-4}
$y_3y_4$ is full.
\end{proposition}
\noindent
{\em Proof.} Suppose $y_3y_4$ is half. Recall $\{x_1,x_2,y_1,y_2\}$ forms a red $K_4$, in which there exists either a red half  $K_3$ or two independent red half  edges from Fact \ref{center-fact} (iv). The former case is impossible from Fact \ref{t=10-fact} (i). For the latter case, since $x_1x_2$ is full, there are another two independent red half  edges in $\{x_1,x_2,y_1,y_2\}$ besides $y_1y_2$ and $y_3y_4$,  contradicting Fact \ref{t=10-fact} (ii).
\hfill$\Box$

\smallskip
By a similar argument as Proposition \ref{t=10-3} by applying Lemma \ref{1-sys}, we  obtain that
\begin{align}\label{y_3,y_4}
N_B(y_3,X)=\{x_2,x_3\}\;\; \text {and}\;\; N_B(y_4,X)=\{x_1,x_5\}.
\end{align}

Suppose first $\{x_1,x_2,y_1,y_2\}$ contains a red half  $K_3$. Since $y_1y_2$ is half and $x_1x_2$ is full, by symmetry, suppose $\{x_1,y_1,y_2\}$ forms a red half  $K_3$. Then $x_2y_4$ and $x_5y_3$ are red full edges from (\ref{y_3,y_4}) and Fact \ref{t=10-fact} (i).
Together with Fact \ref{t=10-inside fact} and (\ref{t=10,deg-x_1x_2}) yield $x_3\in N_R(y_1,X)$ or $x_5\in N_R(y_1,X)$.

If $x_3\in N_R(y_1,X)$, then $x_3y_4$ must be a red half  edge from (\ref{y_3,y_4}); otherwise, $\{x_2,x_3,y_1,y_4\}$ would be a red $K_4$ with a full  $K_3$ on $\{x_2,x_3,y_4\}$, a contradiction. However, then there is a red half  $K_3$ on $\{x_1,y_1,y_2\}$ which is disjoint from the half edge $x_3y_4$, contradicting Fact \ref{t=10-fact} (i).

If $x_5\in N_R(y_1,X)$, then since $x_5y_3$ is a red full  edge, $x_1y_3$ must be a red half  edge from (\ref{y_3,y_4}); otherwise, $\{x_1,x_5,y_1,y_3\}$ induces a red $K_4$ with a full  $K_3$ on $\{x_1,x_5,y_3\}$, a contradiction. Then, $y_2y_3$ would be a red full  edge; otherwise, $x_1y_3y_2y_1x_1$ would form a red half  $C_4$, contradicting Fact \ref{t=10-fact} (iii) since
$x_1y_2$ is half. This implies that $y_2y_4$ is half; otherwise, $\{y_1,y_2,y_3,y_4\}$ would be a red $K_4$ with a red full  $K_3$ on $\{y_2,y_3,y_4\}$, a contradiction. So $y_1y_4$ is full; otherwise, $x_1y_1y_4y_2x_1$ forms a red half  $C_4$, contradicting Fact \ref{t=10-fact} (iii) since $y_1y_2$ is a half  edge. Thus $y_1y_3$ is half; otherwise, $\{y_1,y_2,y_3,y_4\}$ forms a red $K_4$ with a red full  $K_3$ on $\{y_1,y_3,y_4\}$, a contradiction. However, then $\{x_1,y_1,y_3\}$ forms a red half  $K_3$ disjoint from the half  edge $y_2y_4$, contradicting Fact \ref{t=10-fact} (i).

Now we suppose that $\{x_1,x_2,y_1,y_2\}$ contains two independent red half  edges, say $x_1y_1$ and $x_2y_2$ by symmetry. Then we may assume that $x_1y_2$ and $x_2y_1$ are full; otherwise, there exists a red half  $K_3$ in $\{x_1,x_2,y_1,y_2\}$ and we are done. Since $\{y_1,y_2,y_3,y_4\}$ forms a red $K_4$, which contains either a red half  $K_3$ or two independent half  edges by Fact \ref{center-fact} (iv).

First, suppose $\{y_1,y_2,y_3,y_4\}$ contains a red half  $K_3$. Since $y_3y_4$ is full, we may assume that $\{y_1,y_2,y_3\}$ or $\{y_1,y_2,y_4\}$ forms a red half  $K_3$.
If $\{y_1,y_2,y_3\}$ forms a red half  $K_3$, then $x_2y_4,x_3y_4$ are red full  edges from Fact \ref{t=10-fact} (i). It follows from Fact \ref{t=10-inside fact} and (\ref{t=10,deg-x_1x_2}) that
$x_3\in N_R(y_1,X), \; \text{or} \;\; x_5\in N_R(y_1,X).$ This will lead to a contradiction similar as above. Similarly, it is not possible if $\{y_1,y_2,y_4\}$ forms a red half  $K_3$.

Now suppose $\{y_1,y_2,y_3,y_4\}$ contains two independent half  edges. Recall $y_3y_4$ is full. If $y_1y_4,y_2y_3$ are two independent half  edges, then $x_1y_3,x_2y_4$ are full from (\ref{y_3,y_4}); otherwise, $x_1y_1y_2y_3x_1$ or $x_2y_2y_1y_4x_2$ forms a red half  $C_4$, which contradicts Fact \ref{t=10-fact} (iii) since $x_2y_2$ and $x_1y_1$ are half.
It follows from Fact \ref{t=10-inside fact} and (\ref{t=10,deg-x_1x_2}) that $x_3\in N_R(y_1,X)$ or $x_5\in N_R(y_1,X)$.
This will lead to a contradiction similar as above.
Similarly, it is impossible if $y_1y_3,y_2y_4$ are two half  edges.

Now the case when ${\deg_B}( u ) = 4$ is finished, and so Claim \ref{t=10} is complete.\hfill$\Box$

\subsection{Weak stability property for $RT(n,3,7,o(n))$}\label{stb-3-7}


The main idea of the proof of the weak stability property for $RT(n,3,7,o(n))$ is similar to that for $RT(n,3,6,o(n))$, so we only sketch the proof as follows.

Let $\{G_n\}$ be any sequence of asymptotically extremal graphs of the Ramsey-Tur\'{a}n problem $RT(n,3,7,o(n))$, i.e., $G_n$ is $(K_3,K_7)$-free, $\alpha ( G_n )= o( n )$, and
$e( G_n ) = RT( {n,3,7,o( n )} ) + o( n^2 )$. We will show that $\Delta ( G_n,U(n,4,3) ) = o( {{n^2}} )$.
For sufficiently large $n$, we take $\delta,m,\varepsilon$ and $\mu$ with
$0<\frac{1}{n}\ll\delta\ll\frac{\mu^4}{4m^2}$, and $\frac1m\ll\varepsilon\ll\mu^2\ll 1.$
For a vertex $u$ or an edge $uv$, $\chi(u)$ or $\chi(uv)$ denotes the color that $u$ or $uv$ has received.

Using Lemma \ref{b-indent} we obtain a partition $V(G_n)=\sqcup_{k=1}^2  U_k$ so that for $k\in[2]$ every $Y\subseteq U_k$ with $|Y|>\sqrt{\delta} n$ contains an edge of color $k$. Applying Lemma \ref{reg-le} to $G_n$ with $\varepsilon>0$ and $M=\frac{1}{\mu}$ we obtain an integer $M'$ and an $\varepsilon$-regular equitable partition $\sqcup_{i=0}^m V_i$ of $V(G_n)$ with $M\le m\le M'$, which refines the partition $\sqcup_{k=1}^2  U_k$. We may assume $|V_0|=0$ and $|V_i|=\frac{n}{m}$ for $i\in[m]$.
Now, we define a weighted and colored graph $H(\omega_0)$ on $\{v_1,\ldots,v_m\}$ (with a weight function $\omega_0$ and a distribution ${\bf{u_0}}=(u_1,\dots,u_m)$, and allowing also multiple edges), and
(i) for $k\in[2]$, we color a vertex $v_i$ by $\chi(v_i)=k$ iff $V_i\subseteq U_k$;
(ii) for $i\in[m]$, we assign the weight $u_i=\frac{1}{m}$ to the vertex $v_i$;
(iii) for distinct $i,j\in[m]$,
we define the weight of the edge $(v_i, v_j)$ in color $k\in[2]$ as $\omega_{k}(v_i, v_j)$ as (\ref{weight}).
Then
$
 e(G_n)  \le e(H(\omega_0))\left(\frac{n}{m}\right)^2+4\mu {n^2},
$
and so $e(H(\omega_0))\ge \frac{7}{16}m^2-5\mu m^2$.

Let $F_{m,8}$ denote a weighted Tur\'{a}n graph $T_{m,8}$ with $8$ parts $W_1,\ldots, W_8$ such that the weight of each vertex in $F_{m,8}$ is $\frac{1}{m}$, and the weight of each crossing edge in $F_{m,8}$ is $1$.


\ignore{\begin{claim}\label{K_3 nor K_7}
$H(\omega_0)$ contains neither a blue generalized $K_3$ nor a red generalized $K_7$.
\end{claim}}

 We know that $\rho(3,7)={7}/{16}$ from Subsection \ref{v-3-7}, and only $K_8(\omega)$, which is one of the three Ramsey graphs $\{\Gamma_{8}, \Gamma_8', \Gamma_8''\}$ for $r(4,3)$, can achieve the Ramsey-Tur\'{a}n density $7/16$.


Similar to Claim \ref{max}, we have the following claim.

\begin{claim}\label{max-7}
We have $\frac{7}{16}m^2-5\mu m^2\le e(H(\omega_0))\le \frac{7}{16}m^2.$
Furthermore, $e(H(\omega_0))= \frac{7}{16}m^2$ iff $H(\omega_0)$ isomorphic to $F_{m,8}$.
\end{claim}

From Claim \ref{max-7}, if $G_n$ is an asymptotically extremal graph of the Ramsey-Tur\'{a}n problem $RT(n,3,7,o(n))$, then $\Delta(H(w_0),F_{m,8})=o(m^2).$
Recall $U(n,4,3)$ from Construction \ref{con1},  we have $\Delta ( G_n,U(n,4,3) ) = o( {{n^2}} ),$ and so the weak stability property holds as desired. \hfill$\Box$

\section{Concluding remarks and problems}\label{crp}

In this paper, we obtain two classical Ramsey-Tur\'{a}n densities.
We know the values of $\rho(3,q)$ for $q=3,4,5$ and $\rho(4,4)$ from Erd\H{o}s et al. \cite{B3}, and that of $\rho(3,6)$ and $\rho(3,7)$ from Theorem \ref{zhu} and Theorem \ref{zhu-2}.
The values of $\rho(3,q)$ for $q\ge8$ are wide open.

\begin{conjecture}[\cite{B3,kkl}]
For each $q \ge 5$,
$\rho ( {3,2q - 1} ) = \frac{1}{2}(1 - \frac{1}{r(3,q) - 1})$,
and $\rho ( {3,2q}) = \frac{1}{2}(1 - \frac{1}{r(3,q)}).$
\end{conjecture}




Another problem is to determine $\rho ({p_1}, \ldots ,{p_k},\delta )$, which captures more subtle behaviors of the Ramsey-Tur\'an number. Answering two questions proposed by Bollob\'{a}s and Erd\H{o}s \cite{B1}, Fox, Loh and Zhao \cite{flz} showed that $\rho(4,\delta )=\frac{1}{8}+\Theta(\delta)$. Recently,  L\"uders and Reiher \cite{lr} determined the exact value of $\rho({p},\delta )$ for each $p\ge3$. In particular, they proved that for each $p\ge2$, $\rho({2p-1},\delta )=\frac{1}{2}(\frac{p-2}{p-1}+\delta)$, and $\rho({2p},\delta )=\frac{1}{2}(\frac{3p-5}{3p-2}+\delta-\delta^2)$.
For the two-color case, Erd\H{o}s and S\'{o}s \cite{es2} proved in $1979$ that $\rho(3,3,\delta)=\frac{1}{4}+\Theta(\delta)$. Recently,  Kim, Kim and Liu \cite{kkl} determined the exact value of $\rho(3,3,\delta)=\frac{1}{4}+\frac{\delta}{2}$ for all small $\delta>0$, which confirms a conjecture of Erd\H{o}s and S\'{o}s. They also determined the exact values of $\rho(3,4,\delta)$, $\rho(3,5,\delta)$. In the same paper, they gave a construction showing that $\rho(3,6,\delta)\ge \frac{5}{12}+\frac{\delta}{2}+2\delta^2$ and made the following conjecture.

\begin{conjecture}[Kim et al. \cite{kkl}]\label{kkl-2}
For sufficiently small $\delta>0$,
 $\rho(3,6,\delta)= \frac{5}{12}+\frac{\delta}{2}+2\delta^2.$
\end{conjecture}

Given $d,n$ be integers, denote by $F(n,d)$ an $n$-vertex $d$-regular $K_3$-free graph with independence number $d$. Let $S_\delta  \subseteq ( {0,1} )$ consist of all the rationals $\delta$ for which there exists some $F( {n,d} )$ with $\frac{d}{n} = \delta $. From a result of Brandt \cite{Brandt}, we know that $S_\delta$ is dense in $( {0,\frac{1}{3}} )$.

Assume $8$ divides $n$. Let $G$ be a graph obtained from $T_{n,8}$ by putting a copy of $F(\frac{n}{8},d)$, for some $d\in[\delta n-o(n),\delta n]$, in each partite of $T_{n,8}$. It is easy to see that $\alpha(G)\le \delta n$ and $e(G)=\frac{7}{16}n^2+\frac{\delta}{2}n^2+o(n^2)$. Define an edge-coloring $\phi$ of $G$ as follows: $\phi(e)=2$ for all $e\in\bigcup_{i\in[8]}G[X_i]$; $\phi(X_i,X_j)=2$ iff $|i-j|\in[2]$ for all $i,j\in[8]$; all other edges are of color $1$. Then $\phi$ is a $(K_3,K_7)$-free coloring, which implies that $\rho(3,7,\delta)\ge\frac{7}{16}+\frac{\delta}{2}.$


\begin{conjecture}
For sufficiently small $\delta>0$,
 $\rho(3,7,\delta)= \frac{7}{16}+\frac{\delta}{2}.$
\end{conjecture}

Let us pay some attention to the Ramsey-Tur\'{a}n density for non-complete graphs.
Given a forbidden graph $H$, the Ramsey-Tur\'{a}n density $\rho(H)$ for $H$ is defined similarly. An important open problem is to prove a generalization of Erd\H{o}s-Stone Theorem \cite{Erdo-stone}, i.e., $\rho(H)=\rho(p)$ where $p$ depends only on $H$. Define the Arboricity of $H$, denoted $Arb(H)$, as the minimum integer $p$ such that $V(H)$ can be partitioned into $\lceil\frac{p}{2}\rceil$ sets $V_1,\ldots,V_{\lceil\frac{p}{2}\rceil}$ such that $V_i$ spans a forest for each $1\le i\le{\lceil\frac{p}{2}\rceil}$, and $V_{\lceil\frac{p}{2}\rceil}$ spans an independent set if $p$ is odd.
Erd\H{o}s et al. \cite{E-H-3S} proved that $\rho(H)\le\rho(p)$, and the inequality is sharp for odd $p$. Together with $Arb(K_{m,m,m})=5$ for $m\ge3$, Erd\H{o}s et al. \cite{B4} obtained that
$\rho(K_{m,m,m})=\rho(5)={1}/{4}$
 for $m\ge3$.
However, the case for even $p$ is more harder, even the simplest case when $H=K_{2,2,2}$ is not well understood. Since $Arb(K_{2,2,2})=4$, we have $\rho(K_{2,2,2})\le\rho(4)=\frac{1}{8}$. 
\begin{problem}[\cite{B3,B4}]
  Decide if $\rho(K_{2,2,2})=0$?
\end{problem}

\bigskip
\noindent
{\bf Acknowledgement:} The authors would like to thank Hong Liu for the invaluable comments and suggestions which has improved the presentation of the paper greatly.


\end{spacing}
\end{document}